\documentclass[a4paper,11pt, reqno]{amsart}
\usepackage[margin=2cm]{geometry}

\usepackage{graphicx}
\usepackage{amssymb,amsmath,amsthm}
\usepackage{float,geometry}

\usepackage{natbib}
\usepackage[hidelinks,colorlinks = true,linkcolor={blue!80!black},citecolor={blue!50!black},urlcolor={blue!80!black}]{hyperref}
\usepackage{cleveref} 

\usepackage{multirow}
\usepackage{multicol}
\usepackage{adjustbox}
\usepackage{booktabs}
\usepackage{subcaption}
\usepackage{pgfplots}

\usetikzlibrary{shapes,arrows}
\usepackage{setspace}
\usepackage{soul}

\def\R{\mathbb{R}}
\def\D{\mathrm{D}}

\def\S{\mathcal{S}}
\newcommand{\dsum}{\displaystyle\sum}

\newcommand{\dmax}{\displaystyle\max}
\newcommand{\dprod}{\displaystyle\prod}

\newtheorem{prop}{Proposition}[section]

\newtheorem{rmk}{Remark}[section]

\let\origmaketitle\maketitle
\def\maketitle{
	\begingroup
	\def\uppercasenonmath##1{} 
	\let\MakeUppercase\relax 
	\origmaketitle
	\endgroup
}

\begin{document}
	
	\title[]{\Large Insights into Efficiency and Satisfaction Trade-offs in Facility Location Problems with Regional Preferences}

\author[V. Blanco, R. G\'azquez, \MakeLowercase{and} M. Leal]{
{\large V\'ictor Blanco$^{\dagger}$, Ricardo G\'azquez$^{\dagger}$, and Marina Leal$^{\ddagger}$}\medskip\\
$^\dagger$Institute of Mathematics (IMAG), Universidad de Granada, Spain\\
$^\ddagger$Centro de Investigación Operativa, Universidad Miguel Hernández, Spain\\
\texttt{vblanco@ugr.es}, \texttt{rgazquez@ugr.es}, \texttt{m.leal@umh.es}
}
	

\maketitle 	

\begin{abstract}
This paper studies a practical regional demand continuous multifacility location problems whose main goal is to locate a given number of services and entry points in each region to distribute certain products to the users at minimum transportation cost. Additionally, a minimum satisfaction level is required for the customers in each region. This satisfaction is measured through continuous preference functions that reflect the satisfaction degree of each location in the region. We provide a mathematical optimization-based framework for the problem and derive suitable Mixed Integer Second Order Cone optimization models for some interesting situations: norm-based transportation costs for the services to the entry points, and different families of preference functions. Among these preference functions, we highlight those derived from economic production models and distance-based preferences. We conduct an extensive computational study along two main lines: a computational approach, where we provide optimal solutions for up to $500$ demand regions in the single-facility case and up to $50$ for the $p$-facility case; and a qualitative approach, where we analyze whether the incorporation of preferences is statistically significant compared to the case without preferences.
\end{abstract}
		
	
\section{Introduction}

Facility Location is an active area of research in Operations Research whose aim is to find the best position of one or more services satisfying the demand of a set of users. There is a vast amount of literature on location problems, both analyzing their theoretical properties and practical implications. See the monograph by \citep{LaporteNickelSaldanha-da-Gama:2019} for the recent state-of-the-art in the topic.

The different characteristics that may affect a facility location problem imply the use of a wide set of modeling and solution tools of different nature that ranges from Global Optimization, Graph Theory, Mixed Integer Linear and Non-Linear Programming, Stochastic Programming, Robust Optimization, etc. Perhaps, one of the most prominent features of a facility location problem is the domain where the facilities are to be located. On the one hand, \textit{continuous} location problems allow the services to be located in the complete space (or in a region of it) where the users live \citep[see e.g.,][among many others]{BEP14,BEP16,drezner2022continuous,weber1909standort}. In contrast, in \textit{discrete} location, a given finite set of potential positions is provided, and one has to select from such a set the most appropriate facilities \citep[see e.g.,][]{hakimi1964optimum,hakimi1983locating}. In between, in \textit{network} location the facilities are to be located on the metric space induced by a given network \citep[see e.g.,][]{hansen1987single,puerto2018extensive}. These three different frameworks, although all use mathematical optimization as one of the main tools for their analysis, differ on the specific optimization models that need to be solved in each case. Whereas integer linear optimization is the most popular tool for solving discrete and network location problems, global optimization, and convex optimization are to be used in continuous location. 

In recent years, a few papers have appeared analyzing hybridized continuous-discrete facility location problems \citep{blanco2019ordered,blanco2023multi,blanco2022hub} and network-continuous problems~\citep{espejo2023facility} by combining the existence of both discrete/network and continuous type of facilities. One of the families of problems where this combination has been most prominent lately is the one where \textit{neighborhoods} are incorporated into classical discrete/network location problems \citep{blanco2019ordered,blanco2022hub,espejo2023facility} or other types of network design problems~\citep{blanco2017minimum,disser2014rectilinear,espejo2022minimum}. In these problems the elements of the instances (potential facilities, users, edges, and nodes of a network) are assumed to be embedded in a $d$-dimensional space ($\R^d$) and the goal is to construct an optimal solution to the problem by assuming that the optimal location of the services is allowed to be positioned within a given area in $\R^d$ (the neighborhood). These areas may represent the geographical place where the different installations are allowed to be located or the (locational) uncertainty on the exact position where the elements are assumed to be located. The neighborhoods can be modeled either by convex or nonconvex shapes, and obtaining solutions for these problems is a mathematical challenge of great interest in the community since the difficulties of solving both, discrete and continuous location problems, appear. Nevertheless, the study of regions of users in facility location problems is not that recent. \citet{drezner1980optimal} and \citet{carrizosa1998weber} already analyze the single facility Weber problem~\citep{weber1909standort} in case the set of users is identified with a region endowed with a probability measure, and the expected distance to that set of customers is minimized. 

In this paper, we study a facility location problem with neighborhoods of special interest, and that, as far as we know has not been previously addressed. Specifically, we consider a continuous facility location problem with multiple facilities with neighborhoods, where a given finite set of demand regions of general shapes is provided, and the goal is to find the \emph{optimal} positions of a given number of services that minimizes a distance-based measure to the different regions. The peculiarity of the problem that we address is that a function is provided for each of the demand regions to represent the utility of the users in the regions to be serviced. This utility can represent the preference of the users for the different positions within the neighborhood or the population mass defined in the region. For that, each of the demand regions is endowed with a continuous preference function where the demand is serviced from its assigned facility. Thus, one should decide the \textit{entry points} in the regions, i.e. the place in the region where its users receive the service (i.e., a provider hub, a dropbox location, etc). The satisfaction of the users is assured by a minimum preference threshold.

It is usual in facility location to assume that users' preferences are uniquely based on the distance to the service~\citep{disser2014rectilinear}. However, in case the demand is located in regions, the customers distributed all along the region may prefer to locate the service at more accessible places, close to the center of the region, near shopping or industrial areas, or far from conflicting areas. The use of non-distance-based preferences in discrete facility location has already been considered in the literature. \cite{hanjoul1987facility} provide a framework to maximize the preference ordering in the Simple Plant Location Problem. \citet{camacho2014p} and \citet{casas2017solving} propose a bilevel approach for the $p$-median problem, where a leader agent decides on the location of the plants based on distances, but the followers choose the assignment to the open plants based on a different (conflicting) preference function. The discrete gradual covering location problem with preferences, represented as the probabilities that the customers visit the different types of services, is studied by \citet{kuccukaydin2020gradual}. The incorporation of preferences to other types of optimization-based decision models has also be considered in the literature~\citep[see e.g.][]{schoenwitz2017product,wang2024investigation,yura1994production}. There are other works, that are usually classified into the family \textit{equitable/fair} facility location problem, that instead of minimizing the overall distance-based costs from the users to the facilities, consider different cost functions that are minimized, in case the customers are satisfied with a fair distribution of the travel distances from the users' positions to the facilities \citep[see][among some others]{blanco2023fairness,blanco2024intra,espejo2009comparison,chanta2014minimum}. Nevertheless, the study of preferences defined on demand regions in continuous facility problems has not been previously addressed.

Our approach is motivated by the design of central storehouses and containers of e-commerce companies. In the last years, with the increase of the online shopping, the location of hub lockers where consumers can pick up their online orders in each neighborhood has become a common practice. Hub lockers are self-service kiosks that allow users to pick up their packages at a place and time that is convenient for them. Hence, these e-commerce companies must decide where to locate their central storehouses (where products are stored) and the hub lockers (where products are delivered). When deciding these locations, three main ingredients are to be considered: (1) the location of these lockers on certain regions is flexible enough to consider that the decision is space continuous; and (2) the preferences of the users that use the lockers should be taken in to account; and (3) the transportation cost or the distance from the central storehouse to the delivery points, or  the other way round in case of products return, has a remarkable impact in the service time for the users.


Our contributions can be then summarized as follows:
\begin{itemize}
    \item We propose a new mathematical optimization-based framework for continuous multifacility location with demand regions. The decisions to be made are the locations of the facilities and the entry points in each of the demand regions.
    \item We consider general distance-based cost functions between the facilities and the entry points, such that the problem can be cast as a Mixed Integer Second Order Cone Optimization (MISOCO) problem.
    \item The convenience of the entry points in the demand regions is determined by a preference function that indicates how preferred is each point in the region for the users. A threshold is provided, and a minimum preference value is required to determine those entry points to assure the satisfaction of the users.
    \item We provide a general framework to incorporate the preference requirements to our model. We specify different families of preference functions based on two different paradigms, distance-based functions, and economic production models, showing that in all these cases, the problem can be rewritten as a mixed integer conic optimization problem.
    \item We analyze two different cases that may occur in this problem. First, we consider that the demand regions do not overlap, and then, the regions do not share entry points. Secondly, we analyze the general case where the demand regions may intersect and it may be convenient to collocate the entry points in the same position. The latter problem requires further study, and results in a different optimization problem that we develop here.
    \item We design a new data-driven methodology that uses geospatial information to \textit{translate} geographical coordinates into socio-economic features, and that will allow practitioners to apply our methodology to real-world situations.
    \item We report the results of an extensive battery of experiments, analyzing not only the computational limitations of our methodology but also the practical implications of the proposed approach.
\end{itemize}

The paper is organized as follows: Section \ref{sec:preliminaries} is an introduction with the preliminaries of the problem. Section \ref{sec:multi} introduces the continuous facility location problem with preferences model. Section \ref{sec:overlapping} is the explanation of the considered extension, when overlapping between regions occurs. Section \ref{sec:prefs} introduces up to five different preference functions. Section \ref{sec:experiments} reports the extensive battery of experiments. Finally, Section \ref{sec:conclusions} draws some conclusions and future research lines of the paper.

\section{Preliminaries} \label{sec:preliminaries}

The problem studied in this paper is a regional-based extension of the classical Continuous Multifacility Location Problem (CMLP)~\citep[see e.g.][]{BEP16,blanco2019ordered,michelot1987localization}. In the classical single-allocation CMLP, we are given a finite set of demand points in a $d$-dimensional space, $\mathcal{A} = \{a_1, \ldots, a_n\}$, and the goal is to find a given number of new facilities, $p$, such that the sum of the distances from demand points to their closest facilities is minimized. The problem can be abstractly formulated as:
$$
\min_{x_1, \ldots, x_p \in \R^d} \sum_{a\in \mathcal{A}} \min_{j=1, \ldots, p} \D(a,x_j)
$$
where $\D$ is a distance measure in $\R^d$. As in other versions of the CMLP, we consider distances $\D$ based on norms. Specifically, we assume that ${\rm D}(a,b) = \|a-b\|$, for all $a, b \in \R^d$, where $\|\cdot\|$ is a $\ell_\tau$-norm, for $\tau\geq 1$, where 
\begin{equation} \label{eq:ltau_norm}
 \|z\|_\tau = \left\{\begin{array}{cl}
 \left(\dsum_{k=1}^d |z_k|^\tau\right)^{\frac{1}{\tau}} & \mbox{if $\tau<\infty$}\\
 \dmax_{k\in \{1, \ldots, d\}} |z_k| & \mbox{if $\tau=\infty$}\end{array}\right.. 
\end{equation}

The single facility version of the problem above ($p=1$)  can be formulated as a second order cone programming (SOCP) problem and hence is solvable in polynomial time~\citep[see e.g.][]{BEP14,blanco2023minimal}. However, the multifacility case ($p>1$) is known to be NP-hard, even for Euclidean distance in the plane~\citep{megiddo1984complexity}. Some approaches have been proposed to solve the problem exploiting the fact that the problem can be cast as a Mixed Integer Second Order Cone Optimization (MISOCO) problem, even for more general objective functions~\citep{BEP16,blanco2019ordered}, allowing to solve small and medium size instances of the problem.


The first difference between the problem that we address here and the CMLP is that instead of demand points, we consider demand regions. The regions might represent districts in a given town, villages around a county, or any geographical area enclosing a population that needs to be served by the facilities. Thus, instead of measuring the distance between the facility and single points, one has to decide what is the most convenient entry point in the region to measure the distance $\D$ with respect to the location of the facility.

We assume that the demand regions are finite Second Order Cone (SOC) representable sets, i.e., the regions are in the form:
 \begin{equation}
\Big\{z \in \R^d: \|R_k z + T_k\|_2 \leq (c_k)^t z + f_k, \text{ for } k=1, \ldots, m\Big\},
 \end{equation}
\noindent for some $R_k \in \R^{m\times d}, T \in \R^{m}$, $c \in \R^{d\times m}$ and $f \in \R^{m}$.

Recall that SOC-representable sets can be modeled by means of SOC constraints \citep[see][]{BEP14,blanco2023minimal}. Our aproach can be also extended to consider union of SOC-representable sets that can be handled using disjunctive programming tools~\citep{grossmann2002review,espejo2022minimum}. A graphical example of these sets can be seen later in Figure \ref{fig:sol_withoutprefs}.

The second novelty that we consider in this paper is the way the decision-maker decides the best position for the entry point selected in each of the regions. We assume, as usual, that the agent installing the logistic system \textit{pays} for the set-up costs of the facilities and the transportation costs from them to the entry points, but the users pay the transportation cost from them to the pick-up point installed in the region. Thus, in our approach, each region is assumed to be endowed with a \textit{preference function} that represents the satisfaction level of the users for locating the entry point in the different points in the neighborhood.  For each demand region, $R$, we are given a preference function $\Phi_R: R \rightarrow [0,1],$ and we impose that a minimum level of satisfaction, $\Phi_R^*$, is required for the facilities.

\section{Regional Continuous Location Problem with Preferences}\label{sec:multi}

Let $\mathcal{S} = \{S_1, \ldots, S_n\} \subseteq \R^d$ be a set of $n$  SOC-representable demand regions. Each region, $S_i$, is also endowed with a weight $\omega_i >0 $,  for $i \in I:=\{1,\ldots,n\}$, that may represent the population mass or the overall demand for each region. As already mentioned in the previous section, for each region $S_i$ is also given a concave preference function $\Phi_i: S_i \rightarrow [0,1]$ to indicate the degree of preference of customers in $S_i$ for each of the points in the region. Let $\Phi_i^*$ be the minimum desired preference for demand region $i$. We denote by $P_i = \{z \in S_i: \Phi_i(z) \geq \Phi^*\}$, the \textit{feasible} demand region for region $S_i$, i.e., the points achieving this minimum preference threshold, for all $i \in I$.

We assume that transportation costs between facilities are measured using a ($\ell_\tau$ or polyhedral) norm-based distance function $\D$. We consider a single-allocation policy, in which each region will be assigned to exactly one of the facilities, the one whose service point is closest to.

The goal of the \textit{Regional Continuous Location Problem with Preferences} (RCLPP, for short) is to find the position of a set of $p$ new facilities in $\R^d$, and to locate regional entry points in each of the demand regions, minimizing the overall transportation points from the facilities to the service points, and such that the minimum preference level is achieved for all demand regions.

\begin{figure}[ht]
    \begin{subfigure}[b]{.49\linewidth}
        \centering 
        \fbox{\includegraphics[width=0.8\linewidth]{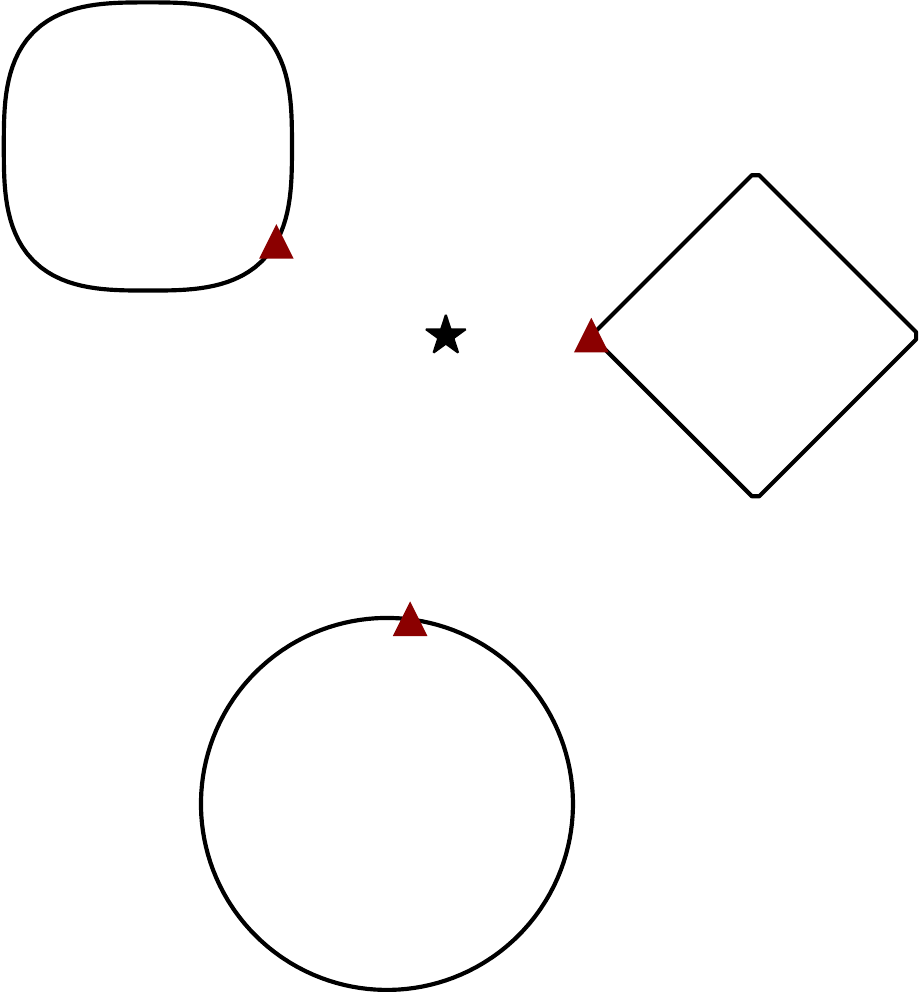}}
        \caption{Without preferences.}\label{fig:sol_withoutprefs}
    \end{subfigure}~\begin{subfigure}[b]{.49\linewidth}
        \centering 
        \fbox{\includegraphics[width=0.8\linewidth]{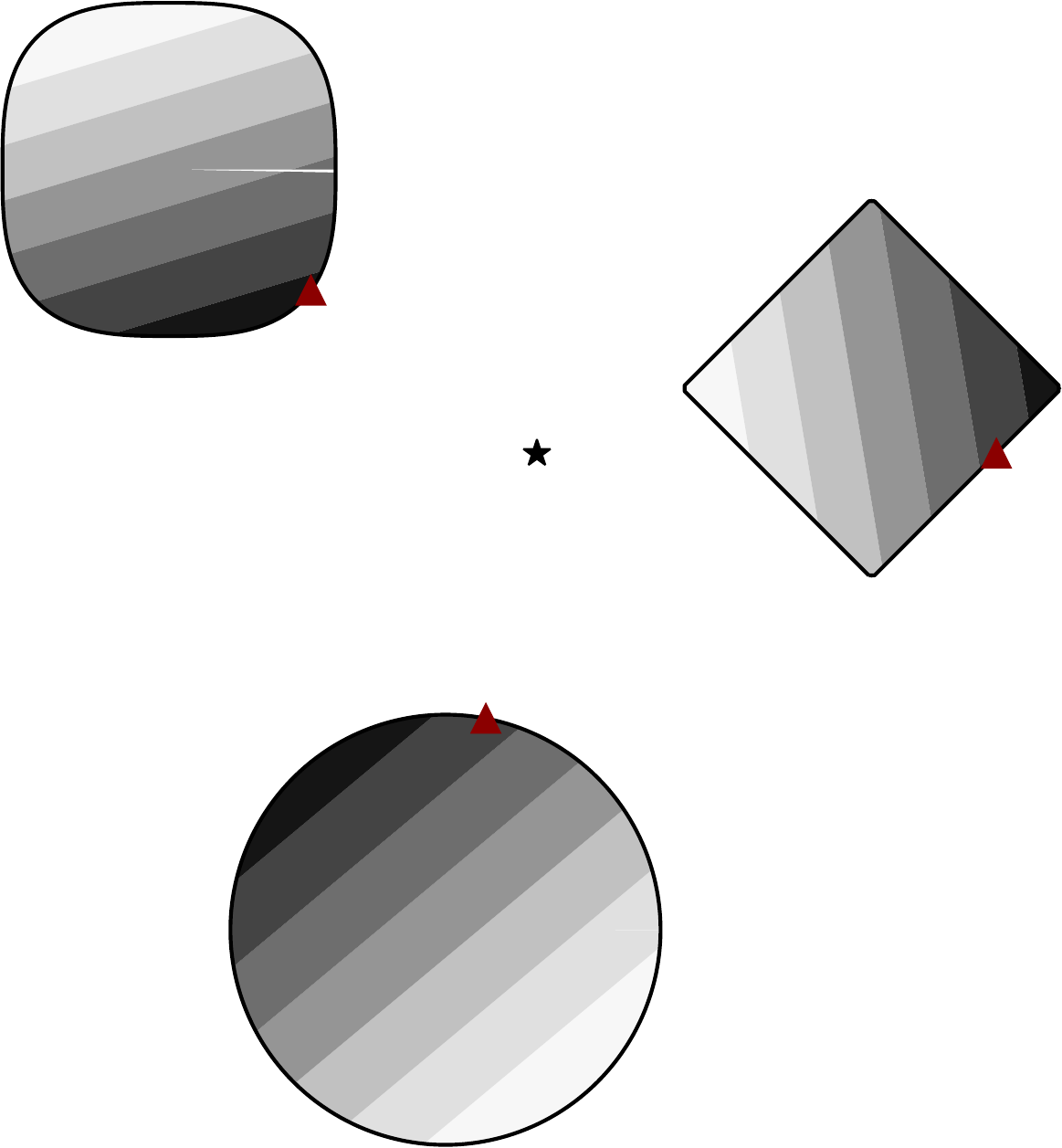}}
        \caption{With preferences.}\label{fig:sol_withprefs}
    \end{subfigure}~
    \caption{Solution for the single-facility case with three SOC-representable regions using different norms. On the left, the solution without preferences is shown, while on the right, the solution with preferences is depicted. The dark red triangles represent the entry points to the regions, and the black star indicates the facility's location, minimizing the distance between itself and the entry points.}\label{fig:solution}
\end{figure}

In Figure \ref{fig:solution} we plot the solutions for an instance of the problem under consideration for a single planar facility and three SOC-representable regions ($\ell_1, \ell_2$, and $\ell_3$-norm-based balls). In the left plot, we draw a solution where no preferences are considered (or the preference threshold is zero), and then the only criterion is the minimization of the transportation costs (Euclidean distance). The point highlighted with a star is the optimal position for the service and the entry points are the red triangles. The only requirement in this case is that the entry points belong to their region. In contrast, in the right plot, we show the case where a preference function and a minimum preference level ($\Phi^*$) are considered for each region.  The contour plot for each preference function within each region is drawn (in this case, it is a linear level preference function, where the darker the color, the better the preference). Although the preference levels are achieved for the left and bottom regions at the same entry points already found for the case where no preferences are considered, the right region requires positioning the point in a different position.

We denote by $J=\{1, \ldots, p\}$ the index set for the facilities to be located.

Note that the problem can be mathematically formulated as:
$$
    \min_{\substack{X \subset \R^d \\ |X|=p \\ a_1 \in P_1, \ldots, a_n \in P_n}} \sum_{i \in I} \min_{x \in X} \D(a_i, x)
$$

In what follows, we derive a suitable mathematical optimization formulation to be solved with off-the-shelf software, and then, being used for practitioners.

In Table \ref{table:vars} we detail the decision variables that we use in our model. 

\begin{table}[ht]
\centering
\begin{tabular}{|rl|}\hline
Variables & Decision\\\hline
$x_j \in \R^d$ & Coordinates of the $j$th facility, $j \in J$.\\
$a_i \in \R^d$ & Coordinates of the service point in the $i$th region, $i\in I$.\\
$y_{ij} \in \{0,1\}$ & $1$ if region $i$ is allocated to facility $j$, and $0$, otherwise, $i\in I$, $j\in J$.\\\hline
\end{tabular}
\caption{Decision variables used in our mathematical optimization model. \label{table:vars}}
\end{table}

With these variables, we propose the following optimization formulation for the problem:
\begin{subequations}\label{model:multi}
\begin{align}
\min \quad & \dsum_{i\in I} \dsum_{ j \in J} \omega_i \D(a_i,x_j) y_{ij} \label{multi:of}\\
\mbox{s.t.}\quad 
&\sum_{ j \in J} y_{ij}=1,\quad \forall i \in I, \label{multi:assigment}\\
& a_i \in P_i,\quad \forall i \in I,\label{multi:dom-a}\\
& x_j \in \R^d,\quad \forall  j \in J, \label{multi:dom-x}\\
& y_{ij}\in \{0,1\},\quad \forall i \in I, j \in J. \label{multi:dom-y}
\end{align}
\end{subequations}

In this formulation, the objective function represents the minimization of the transportation cost between the facilities and the entry point if the entry point $i$ is assigned to facility $j$ ($y_{ij} = 1$).  The set of constraints  \eqref{multi:assigment} ensures that each region is assigned to exactly one facility. Constraints \eqref{multi:dom-a}  forces the location of each regional center to both belong to the demand region and to achieve the preference threshold required by this region.

Note that \eqref{model:multi} is a discrete and non-linear model. The nonlinearities of the model appear both in the objective function \eqref{multi:of} and in constraint \eqref{multi:dom-a}. The objective functions is the sum of the products of nonlinear functions (the distances from facilities) by binary variables ($y$). In this case, one can linearize the expression by replacing each term $\D(a_i,x_j)y_{ij}$ by a new auxiliary continuous variable $d_{ij}$. To assure that the $d$-variables attain the desired values, the following family of constraints needs to be added to the model: 
\begin{equation}\label{multi:linearization}
d_{ij} \geq \D(x_{j}, a_{i}) - \Delta^1_i (1-y_{ij}), \quad \forall i \in I, j \in J,
\end{equation}
\noindent where, for each $i\in I$, $\Delta^{1}_i > \D(q_i,q_{i'})$, for all  $q_i \in S_i, q_{i'} \in S_{i'}$, for $i, i' \in I$. The model still remains non-linear due to the non-linearity of the norms. Note that the sets $P_1, \ldots, P_n$ are defined as the intersection of a SOC-representable convex sets and the desired preference region based on the given threshold. Thus, in general, these sets cannot be described by linear sets of inequalities, implying another degree of nonlinearity in the formulation above. However, since the sets $S_i$ are SOC-representable, we construct in the next section, preference functions that preserve the SOC-representability of the demand regions, by intersecting them with SOC-representable sets induced by the preference functions. With these assumptions, the model would allow a SOC representation, resulting in a Mixed Integer Second Order Cone (MISOCO) programming model which can be solved by off-the-shelf optimization solvers \citep[see][]{BEP14,blanco2023minimal}. 

Observe that if a single overall preference function is provided instead of one for each demand region (that is, a general overall satisfaction for all the users), one may replace the threshold constraints in our model by $ \Phi(a_1, \ldots, a_n) \geq \Phi^*$ for any given preference aspiration level $\Phi^*$. The overall preference function $\Phi: \otimes_{S_i \in \S} S_i  \rightarrow [0,1]$ can be, for instance, the minimum overall preference functions for the different regions, the mean of all preferences, or others. 

The single facility problem is a particular case of our problem that can be solved more efficiently, since the resulting model is continuous, as detailed in the following result.
\begin{rmk}
The single-facility model can be formulated as a continuous mathematical optimization problem with SOC constraints. In this case, the binary allocation variables are no longer required, and only the continuous remain in the resulting model:

 \begin{subequations} \label{P1}
\begin{align}
\min_{\mathbf{a},x}\quad & \dsum_{i\in I} \omega_i \D(a_i,x)\label{single:of}\\
\mbox{s.t.} \quad & a_i \in P_i,\quad \forall i\in I, \label{single:dom-a}\\
& x \in \R^d. \label{single:dom-x}
\end{align}
 \end{subequations}

It corresponds to a continuous model with SOC constraints in distance constraints and $\S$-regions.
\end{rmk}

In the general location problem that we address here, the solution space for the variables $\mathbf{x}$ and $\mathbf{a}$ is a subset of a continuous $d$-dimensional space. Thus, $\mathbf{x}$ can take any value in $\R^d$ and $\mathbf{a}$ any point in $\S \subset \R^d$ that satisfies the minimum required preference level. This case generalizes the possible combinations when we consider as solution spaces a discrete space, either for $\mathbf{x}$, for $\mathbf{a}$, or both. The following remarks study the cases where the solution space is/are discrete.
\begin{rmk}
The discrete version of the problem under analysis can be considered by assuming that the servers and/or the entry points are selected from a finite set of points. In what follows we describe these situations:
\begin{itemize}
  \item If the solution space for $\mathbf{a}$ is discrete, this results in a finite set of possible locations for the entry points within the regions. We can address this case using binary variables $t_{ik}$ taking the value 1 if in region $i \in I$ we open the entry point at $k$, with $k \in \{1, \ldots K_i\}$, where $K_i$ is the total number of entry points candidates in region $i$. The model can still be formulated as a MISOCO programming model.
\item 
 If the solution space for $\mathbf{x}$ is discrete, that means that we have a finite set of potential locations for the facilities. The binary variable $r_{jk}$ can be included to adapt the model to this particular situation, taking the value of these variables $1$ if the facility $j$ is located at a potential location $k$, with $j \in J$ and $k \in \{1,\ldots,K\}$, which is $K$ the number of potential locations for the facilities. It results in a MISOCO programming model as before.
\item If both spaces are discrete, both described binary variables, $\mathbf{t}$ and $\mathbf{r}$,  need to be included. SOC constraints are removed since the distances can be computed beforehand. Resulting in this case in a full integer programming (IP) model.
\end{itemize}
\end{rmk}

\section{Regional Continuous Location Problem with Preferences and Collocation}\label{sec:overlapping}

Different demand regions may overlap, i.e., share nonempty spaces where their service points can be located in the same position (\emph{collocation}). Although the model presented above does not consider this situation, it may be advisable since transportation costs can be reduced if the preference thresholds are satisfied. In this section, we provide a modification of our approach that considers this situation and allows to find optimal solutions under this framework. Observe that with the model presented in the previous section collocating entry points is not recommendable since the distance from the service to the entry point would be duplicated. 

Let us consider the following instance for the problem: Let $P_1, \ldots, P_n \subseteq \R^d$ the preference-constrained demand regions such that $\bigcap_{i=1}^n P_i \neq \emptyset$. Let us assume that $a_i \in P_i$ is the service point of the $i$th demand region. In case $a:=a_1= \cdots = a_n \in \bigcap_{i=1}^n P_i$, then from the previously presented model the transportation cost from the facility where they are allocated, say $x_j$, is $n \D(a, x_j)$. Although this is an extreme case where all the regions intersect, it may happen that two regions intersect, and the same situation occurs.

\begin{figure}[ht]
    \begin{subfigure}[b]{.49\linewidth}
        \centering 
        \fbox{\includegraphics[width=0.8\linewidth]{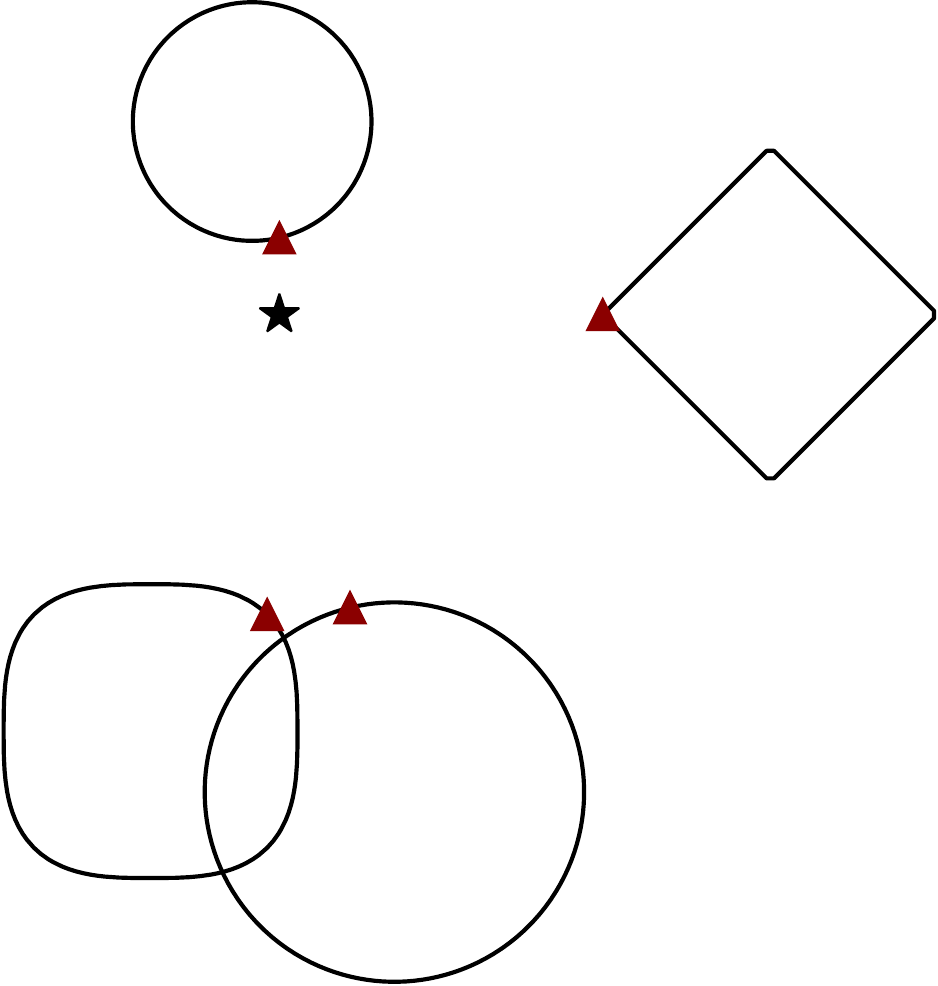}}
        \caption{Without collocation.}\label{fig:without_overlapping}
    \end{subfigure}~\begin{subfigure}[b]{.49\linewidth}
        \centering 
        \fbox{\includegraphics[width=0.8\linewidth]{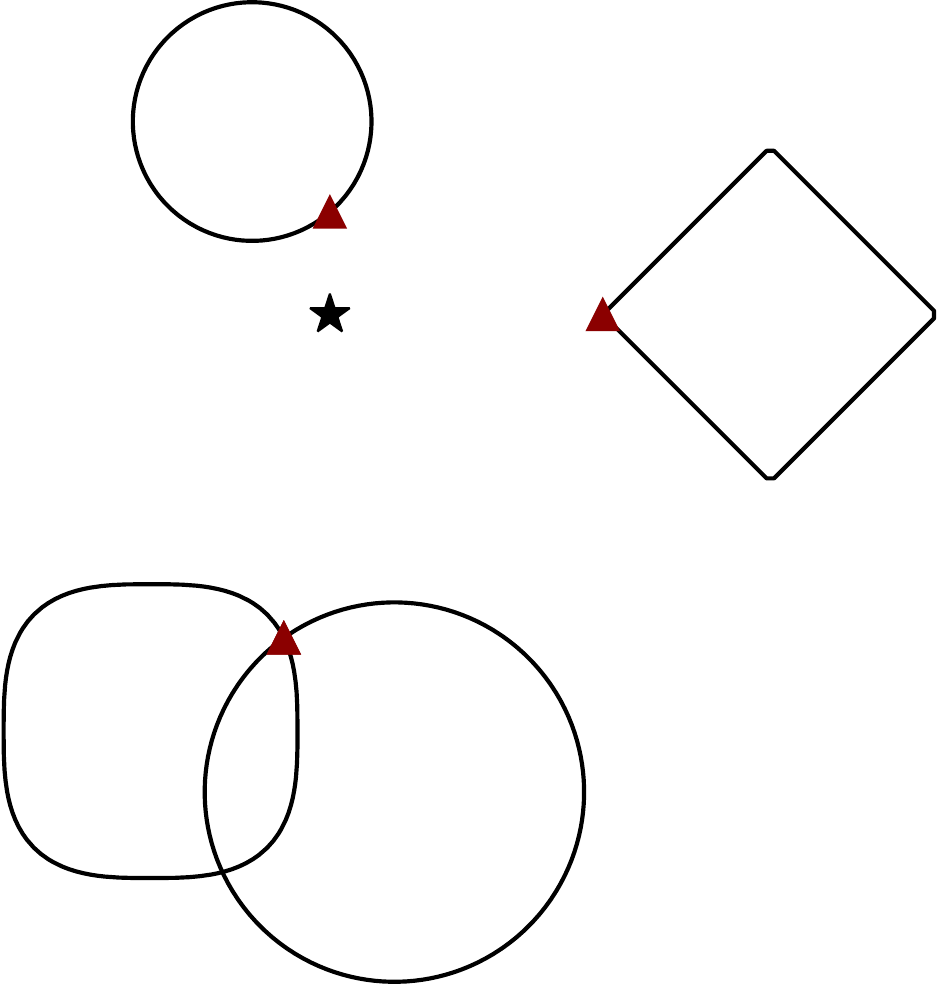}}
        \caption{With collocation.}\label{fig:with_overlapping}
    \end{subfigure}~
    \caption{Solutions for the single-facility case with four regions using different norms without preference functions. On the left, the solution without considering collocation, while on the right, the solution with collocation is depicted. Dark red triangles represent the entry points to the regions, and the black star indicates the facility's location, minimizing the distance between itself and the entry points.}\label{fig:overlapping}
\end{figure}

To show the impact of overlapping regions, in Figure \ref{fig:overlapping} we present solutions for an instance for a single-facility problem with four regions (and no preferences). In the left plot, we draw the solution for the case where collocation is not contemplated in the procedure. Thus, one entry point is positioned for each region. However, the two bottom regions intersect, and it may be convenient to collocate their entry points to avoid almost doubling the transportation costs from the service to these two points. In the right plot, we draw the solution in case collocation is considered. In this case, instead of four entry points, the solution constructs three different positions for the entry points, collocating those for the two bottom regions. This consideration has a direct impact on the objective value, in this case, a reduction of around 30\%.

To overcome this situation we propose an alternative model that allows for the collocation of service points (in case the regions intersect) by checking the pairwise coincidences of the service points in $P_1, \ldots, P_n$, and only account for the distance from the assigned service to one of the entry points. Thus, in the new formulation we need to incorporate decision variables to detect whether the entry points of different intersecting regions coincide, and in such a case, only account in the objective function for one of the transportation costs from the corresponding facility. Then, additionally, to the variables already used in the previous model, the following set of binary variables are considered:
$$
z_{ii'} = 
\begin{cases}
1 & \text{if } a_i = a_{i'}, \\
0 & \text{otherwise}
\end{cases}
\quad \text{for all } i, i' \in I \text{ with } P_i \cap P_{i'}
$$

With the above set of variables, we define the \textit{actual} distance from a facility located at $x_j$ and service point $a_i$ as:
\begin{equation}\label{eq:d}
 d_{ij} = \D(a_i, x_j) y_{ij} \prod_{\substack{i'\in I: \\ i'< i}} (1-z_{ii'}).
\end{equation}
\noindent That is, in case $a_i$ is allocated to facility $x_j$ ($y_{ij}=1$) and there is no $i'<i$ with the same location for the service point in $P_{i'}$ that $a_i$, then, the distance is $\D(a_i,x_j)$. Otherwise, in case $a_i$ is not allocated to $x_j$ or there exists $i'<i$ with $a_i=a_{i'}$, the distance is accounted as $0$ in the objective function.  By convention, in case of collocation, the distance to the assigned facility is calculated using the smallest-index collocated entry point. Thus, $\dprod_{\substack{i'\in I: \\ i'< i}} (1-z_{ii'})$,  
takes the value 1 if and only if $i$ is the smallest index where the entry points coincide.

The adequate definition of the $z$-variables above is assured by the following constraints:
\begin{equation}
    \|a_i - a_{i'}\|_1 z_{ii'} \leq 0 , \quad \forall i,i' \in I, i' < i \mbox{ with $P_i\cap P_{i'} \neq \emptyset$}.
\end{equation}
\noindent That is, if $z_{ii'}=1$, the points $a_i$ and $a_{i'}$ are forced to coincide. Otherwise, the constraint is redundant. Note that the use of the $\ell_1$-norm allows the rewriting through a set of linear constraints.

With the above considerations, we propose the following mathematical optimization formulation for the problem that we name the \textit{Collocation Regional Continuous Location Problem with Preferences} (C-RCLPP), for short:

\begin{subequations}
\begin{align}
\min \quad & \sum_{i\in I}\sum_{j\in J} \D(a_i, x_j) y_{ij} \prod_{\substack{i'\in I: \\ i'< i}} (1-z_{ii'})\label{overlapping1:of}\\
\mbox{s.t.}\quad 
&\sum_{ j \in J} y_{ij}=1,\quad \forall i \in I, \label{overlapping1:assigment}\\
& \| a_{i} - a_{i'} \|_{1} z_{ii'} \leq 0, \quad \forall i,i' \in I, i' < i \mbox{ with $P_i\cap P_{i'} \neq \emptyset$}, \label{overlapping1:entrance}\\
& a_i \in P_i,\quad \forall i \in I,\label{overlapping1:dom-a}\\
& x_j \in \R^d,\quad \forall  j \in J, \label{overlapping1:dom-x}\\
& y_{ij}\in \{0,1\},\quad \forall i \in I, j \in J, \label{overlapping1:dom-y}\\
& z_{ii'}\in \{0,1\},\quad \forall i \in I, i' \in I. \label{overlapping1:dom-z}
\end{align}
\end{subequations}

Observe that the above formulation is, again, a mixed integer nonlinear optimization problem, both because the objective function (which is a polynomial function) and constraints \eqref{overlapping1:entrance}. In the following result we prove that the formulation can be adequately rewritten to be solved using off-the-shelf MISOCO solvers.

\begin{prop}
 The above mathematical optimization is equivalent to the following MISOCO problem:
\begin{subequations}
\begin{align}
 \min \quad & \sum_{i\in I}\sum_{j\in J} d_{ij}\\
\mbox{s.t.}\quad 
&\sum_{j \in J} y_{ij}=1,\quad \forall i \in I, \label{overlapping2:assigment}\\
& d_{ij} \geq \D(a_{i}, x_j) - \Delta^{1}_i(1-y_{ij}+\dsum_{\substack{i'\in I: \\ i'< i}} z_{ii'}), \quad \forall i \in I, j \in J, \label{overlapping2:linearization}\\
& \| a_{i} - a_{i'} \|_{1} \leq \Delta_{ii'}^{2}(1-z_{ii'}), \quad \forall i,i' \in I, i' < i \mbox{ with $P_i\cap P_{i'} \neq \emptyset$}, \label{overlapping2:entrance}\\
& a_i \in P_i,\quad \forall i \in I,\label{overlapping2:dom-a}\\
& x_j \in \R^d,\quad \forall  j \in J, \label{overlapping2:dom-x}\\
& y_{ij}\in \{0,1\},\quad \forall i \in I, j \in J, \label{overlapping2:dom-y}\\
& z_{ii'}\in \{0,1\},\quad \forall i \in I, i' \in I, \label{overlapping2:dom-z}
\end{align} 
\end{subequations}
\noindent where, for each $i\in I$, $\Delta^{1}_i > \D(q_i,q_{i'})$, for all  $q_i \in S_i, q_{i'} \in S_{i'}$, for $i, i' \in I$, and $\Delta_{ii'}^2 = \max_{q_i \in P_i, q_{i'}\in P{i'}} \{\|q_i-q_{i'}\|_1\}$ for $i, i' \in I$ with $P_i\cap P_{i'} \neq \emptyset$.
\end{prop}
\begin{proof}
The adequate definition of each of the distances in the objective function, where $d_{ij}$ is defined as in \eqref{eq:d} is ensured by constraint \eqref{overlapping2:linearization}. Specifically, the expression in the right-hand side is nonnegative just in case $y_{ij}=1$ ($a_i$ is allocated to $x_j$) and $\dsum_{\substack{i'\in I: \\ i'< i}} z_{ii'} = 0$ (that is, in case there is no smaller index to $i$ for which $z_{ii'}=1$). Since all the distances are nonnegative and the objective function is minimized, the distance $d_{ij}$ takes then either value $0$ or $\D(a_i,x_j)$.

Analogously, constraint \eqref{overlapping2:entrance} is a linearization of constraint $\| a_{i} - a_{i'} \|_1z_{ii'} \leq 0$.
\end{proof}

\section{Suitable preference functions}\label{sec:prefs}

Given the definition of the set $P_i = \{z \in S_i: \Phi_i(z) \geq \Phi^*\},\; i \in I$, this section is devoted to propose five  different families of preference functions ($\Phi$) that can be incorporated to our models: linear, based on distances to points of interest, and three based on economic production models. 
Moreover, we also show that determining thresholds $\Phi_i^*$ depends on the choice of the different parameters in the preference functions. In order to define a more intuitive threshold for the decision-maker, we show in this section the normalization of each function to take values between 0 and 1, which are handler for interested decision-makers. Finally, we derive suitable constraints representing the preference to be incorporated into the model.

\subsection*{Linear preferences}

Linear preferences are useful in situations in which the separation with respect to a linear border implies a linear increase or decrease of the demand. This function assumes that the preference or mass accumulation of demand increases linearly in the direction of $\gamma_i$ and it is defined as:
\begin{equation}\label{eq:pref_lineal}
 \Phi^L_i(x) =\gamma_i x + \gamma^{0}_{i},\; \forall x \in S_i, 
\end{equation}
\noindent for some $\gamma_i \in \R^d$ and $\gamma_{i}^{0}\in \R$, $\forall i \in I$. 
In this case, the set of constraints \eqref{multi:dom-a} in Problem \eqref{model:multi} can be rewritten as:
$$
\gamma_i^t x_i + \gamma_0 \geq \Phi_i^*, \forall i \in I, x_i \in S_i.
$$

Figure \ref{fig:prefs_L} illustrates an example of the linear preference within the regions. Note how the preference ranges from a lighter to a darker color, indicating the transition from the lowest to the highest minimum preference level for locating the entry point (from 0 to 1 in the normalized scale).

\subsection*{Distance-based preferences}

These preference functions are useful in cases in which the closeness to certain usually frequented points in the region, as for instance supermarkets, schools, shopping centers, petrol stations, etc., is preferred for the users of serving points. Thus, for each demand region $S_i$, let $B_{i1}, \ldots, B_{iK_i} \in \R^d$, $K_i > 0$ reference points whose closeness from the serving points to them is desired. The distance-based preference functions are defined as:
\begin{equation}
 \Phi_i^{D}(x) = - \dsum_{k=1}^{K_i} \lambda_{ik} \D(B_{ik},x),
\end{equation}
\noindent where $0\leq \lambda_{ik} \leq 1$ with $\dsum_{i=1}^{K_i} \lambda_{ik}=1$. The closest $x$ to the reference points $B_{ik}$, the better. Note that the negative sign in the preference function is due to the fact that the preferences are usually maximized, and in this case, the higher the distance, the worse; hence, this negative sign corrects this issue.
Under these distance-based preference function, constraints  \eqref{multi:dom-a} are replaced by: 
$$
- \dsum_{k=1}^{K_i} \lambda_{ik} \D(B_{ik},x_i) \geq \Phi_i^*, \forall i\in I, x_i \in S_i.
$$

\begin{figure}[h!]
    \begin{subfigure}[b]{.49\linewidth}
        \centering 
        \fbox{\includegraphics[width=0.8\linewidth]{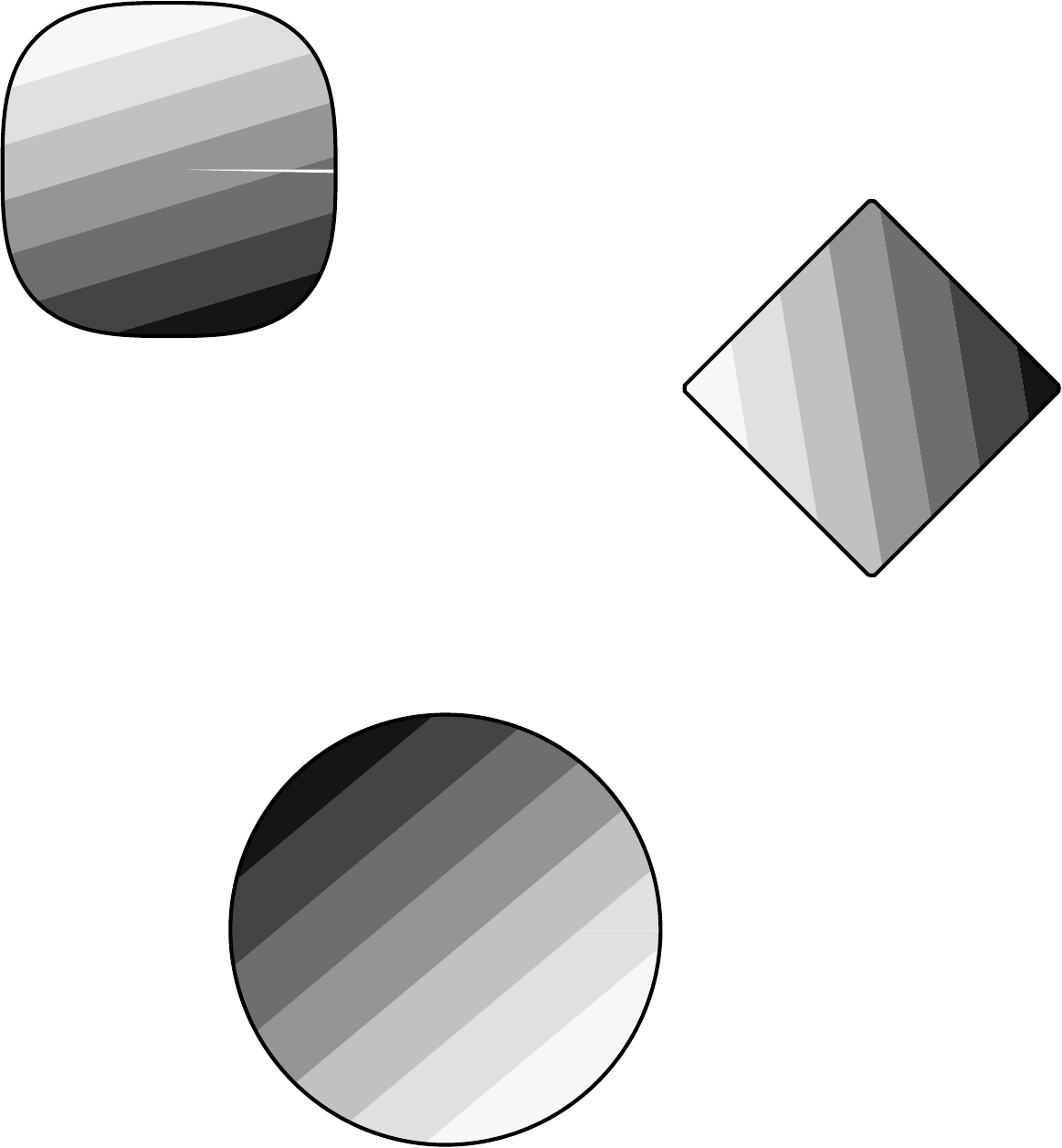}}
        \caption{Linear preference function.}\label{fig:prefs_L}
    \end{subfigure}~\begin{subfigure}[b]{.49\linewidth}
        \centering 
        \fbox{\includegraphics[width=0.8\linewidth]{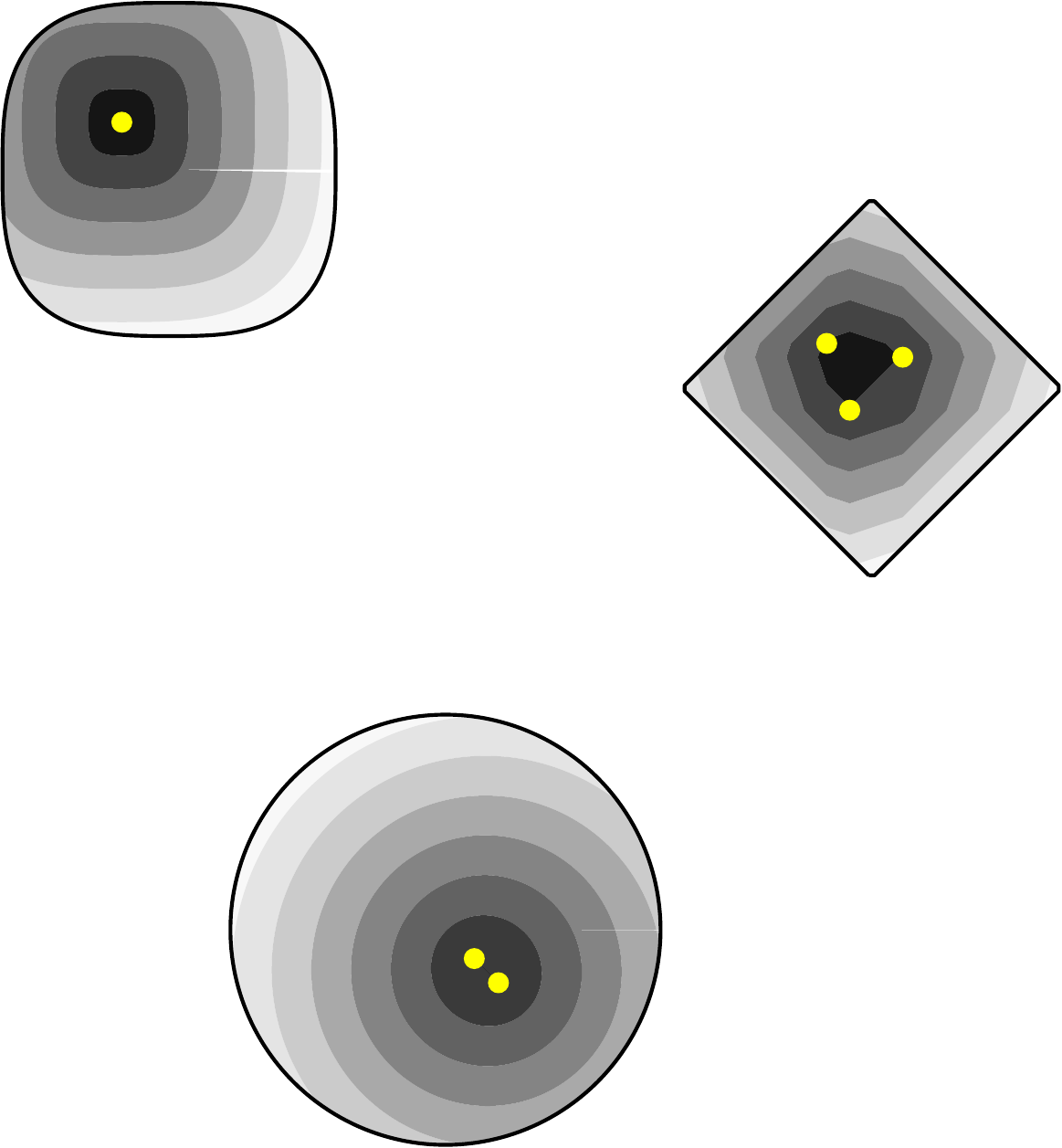}}
        \caption{Distance preference function.}\label{fig:prefs_D}
    \end{subfigure}~
    \caption{Example of two preference functions. From lightest to darkest, from lowest to highest minimum preference level.}\label{fig:prefs_LD}
\end{figure}

Figure \ref{fig:prefs_D} illustrates the distance preference function within the regions. We denote in yellow the reference points in the region. Again, the preference level ranges from a lighter to a darker color, indicating the transition from the lowest to the highest minimum preference level for locating the entry point (from 0 to 1 in the normalized scale). Note that the highest preference level is closer to the yellow points. 

\subsection*{Preferences based on production models}

In what follows we detail a different family of preference functions to represent regional preferences based on classical production models, namely Cobb-Douglas (CD), Constant Elasticity of Substitution (CES), and $m$-factor Leontief (LF). A production model is derived to \textit{estimate} a production output in terms of different inputs. Classical models use as input the available resources (workforce, raw material, equipment, etc.) in order to adequately determine the level of production of certain goods. These models allow to return as output the degree of production as a function of the different inputs. We adapt here these three different types of production models in order to determine the preference of each point of the region given different characteristics of their position. The preference functions that we consider for each region $S_i$ and any location $x\in S_i$ are in the form:
\begin{equation}
 \Phi_i^{PM}(x) = f_i(g_i(x)),
\end{equation}
\noindent where $f_i: \R^m \rightarrow \R_+$ is a (concave) production function and $g_i: S_i\rightarrow \R^m$ is an affine transformation of the coordinates of the points in $S_i$ to different characteristics of the points in $S_i$ that affect to the preference of the users. 

For each region $S_i$ and any $x \in S_i$, we consider the following families of production models: 

\begin{description}
    \item[Cobb-Douglas (CD).] The CD model introduced by \citet{cobb1928theory}, has been widely analyzed in Economics when studying the relationship between the production ($P$), the labor ($L$) and the capital ($K$) when producing certain good, with a function if the form $P = \kappa L^\beta K^{1-\beta}$. Where the degree $\beta$ (also known as elasticity) allows to model the production response for changes in the levels of either labor or capital ($1\%$ increase of labor implies a $\beta \%$ increase in production and $1\%$ increase in capital results in a $(1-\beta)\%$ in production). It is defined as follows,
    \begin{equation}\label{eq:CD}
     \Phi_i^{CD}(x) = \prod_{j=1}^m g_{ij}(x)^{\beta_{ij}},
    \end{equation}
    \noindent for $\beta_1, \ldots, \beta_{m} \geq 0$ with $\dsum_{j=1}^{m} \beta_{ij}=1$. In this case, the product in the right-hand side represents how satisfying is the location $g(x)$. An increasing of $1\%$ in the $j$th feature implies a global satisfaction increasing of $\beta_{ij}\%$.
    
    \item[Constant Elasticity of Substitution (CES).]
    The CES utility function was introduced by \citet{solow} to model production in terms of labour and invested capital. It is defined as,
    \begin{equation}\label{eq:CES}
     \Phi_i^{CES}(x)= \left(\dsum_{j=1}^{m} \beta_{ij} g_{ij}(x)^{\tau_i}\right)^{\frac{1}{\tau_i}},
    \end{equation}
    \noindent for $\beta_1, \ldots, \beta_{m} \geq 0$ with $\dsum_{j=1}^{m} \beta_{ij}=1$, $0<\tau_i\leq 1$. These $\beta$-parameters determine the degree of substitutability between inputs, and they are also known in the literature as elasticity of substitution.
    
    \item[$m$-factor Leontief (LF).]  The Leontief function given by \citet{leontief1941}, and the extension to $m$ factors is,
    \begin{equation}\label{eq:GL}
     \Phi_i^{L}(x) = \min\left\{\frac{g_{i1}(x)}{\beta_{i1}}, \ldots, \frac{g_{im}(x)}{\beta_{im}}\right\} 
    \end{equation}
    \noindent where $\beta_{ij} > 0$ indicates the (minimal) amount of factor $j$ required to produce a single unit of preference for the $i$th-region.

\end{description}

The functions described above have been widely used to measure the effect of different features on the production of a good. Production is nothing but the result of combining different types of inputs in order to make something for consumption. In our framework, preference of a position can be seen as some type of \textit{production} that is affected by different attributes related to the regional situation, as for instance the closeness to a bus/railway station, the prices of houses, the number of shops around the point, etc. We measure these features by means of the $g$-functions above. 

In what follows we describe a data-driven methodology that uses geospatial information to \textit{translate} geographical coordinates into socio-economic features that may affect the preference functions considered in this work.

Let us consider a feature of interest to be measured in demand region $S_i$. Our approach consists of subdividing the region in $L$ regions based on the values of this feature. For the sake of simplicity and applicability, we consider piecewise constant functions in the form:
\begin{equation}
 g_{i}(x) = \begin{cases}
\delta_{i1} & \mbox{if  $x\in s_{i1}$},\\
\delta_{i2} & \mbox{if  $x\in s_{i2}$},\\
\vdots & \vdots\\
\delta_{iL} & \mbox{if  $x\in s_{iL}$}.
\end{cases}
\end{equation}
where $s_{i1}, \ldots, s_{iL} \subseteq S_i$. Each of these regions is then assigned a value for the feature of interest. We assume that $\bigcup_{l=1}^L s_{il} = S_i$. This simple modeling approach allows to indicate the zones of the region with a larger number of parking spaces, or those with more bus stops, malls, etc. 

We assume that the subregions are in the form:
$$
s_{il} = S_i \cap \{z\in \R^d: A_{il}z \leq b_{il}\}
$$
with $A_{il} \in \R^{m_{il}\times d}$, $b\in \R^{m_{il}}$. That is, the $s$-regions is a linear subdivision of the region.

This function can be incorporated to our mathematical optimization formulations via the following linear constraints and auxiliary binary variables:
\begin{align}
 & g(a_i) = \dsum_{l=1}^{L_i} \delta_{il} \xi_{il}, \; \forall i \in I, \label{eq:subdiv1} \\
 & A_{il} a_i \leq b_{il} + \Delta^{3}_i (1-\xi_{il}), \; \forall i \in I,\;,\; l \in \{1,\ldots,L_i\}, \label{eq:subdiv2} \\
 & \dsum_{l=1}^{L_i} \xi_{il} = 1,\; \forall i \in I,\; \label{eq:subdiv3} \\
 & a_i \in S_i, \; \forall i \in I, \label{eq:subdiv4} 
\end{align}
\noindent where $\Delta^{3}_i > \max \{A_{il}a - b_{il} : a \in S_i\}$ for all $i\in I$.

The first set of constraints \eqref{eq:subdiv1} allows one to determine the values of the $g$-functions, constraints \eqref{eq:subdiv2} fix in which subdivision is the variable $a_i$, the family of constraints \eqref{eq:subdiv3} says that $a_i$ have to belong to only one subdivision, and finally constraints \eqref{eq:subdiv4} indicate the belonging to the region $S_i$.

In this way, one can discretely represent the values of a single feature of interest. It is usual to assess the preference of a certain zone in the region via different features, i.e., combining the availability of parking slots, bus stops, number of ATMs, etc to get an overall satisfaction value for each point in the demand region. Thus, one can construct as many of these functions as desired, as different layers with the same underlying demand region. In that case, each of these functions must be represented as the single $g$-function above (with a set of binary variables for each of them). Then, all of them are combined using the production models above to construct preference functions from linear partitions of the regions.

\begin{figure}[ht]
    \begin{subfigure}[b]{.33\linewidth}
        \centering 
        \fbox{\includegraphics[width=0.8\linewidth]{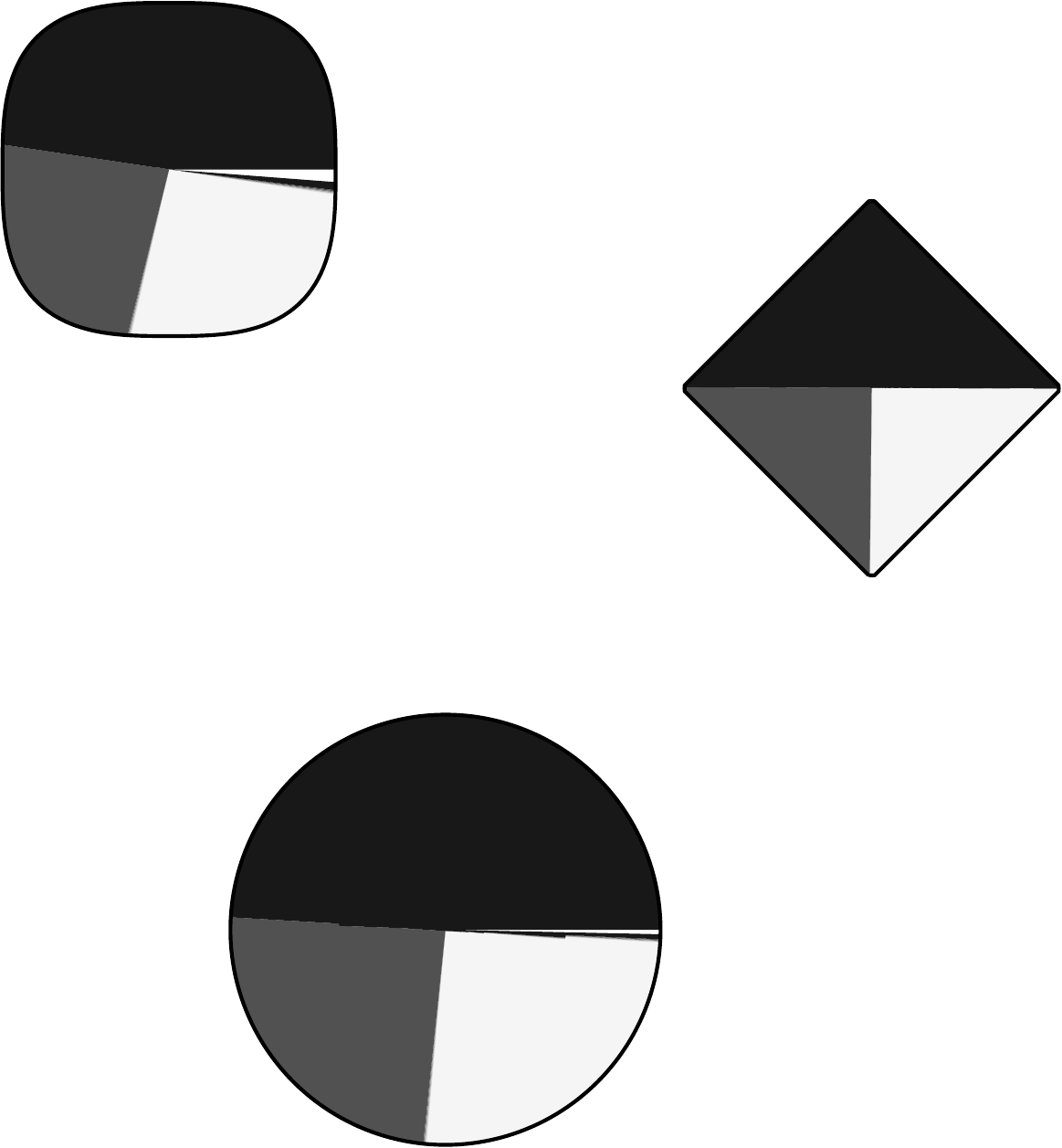}}
        \caption{CD preference function.}\label{fig:prefs_CD}
    \end{subfigure}~\begin{subfigure}[b]{.33\linewidth}
        \centering 
        \fbox{\includegraphics[width=0.8\linewidth]{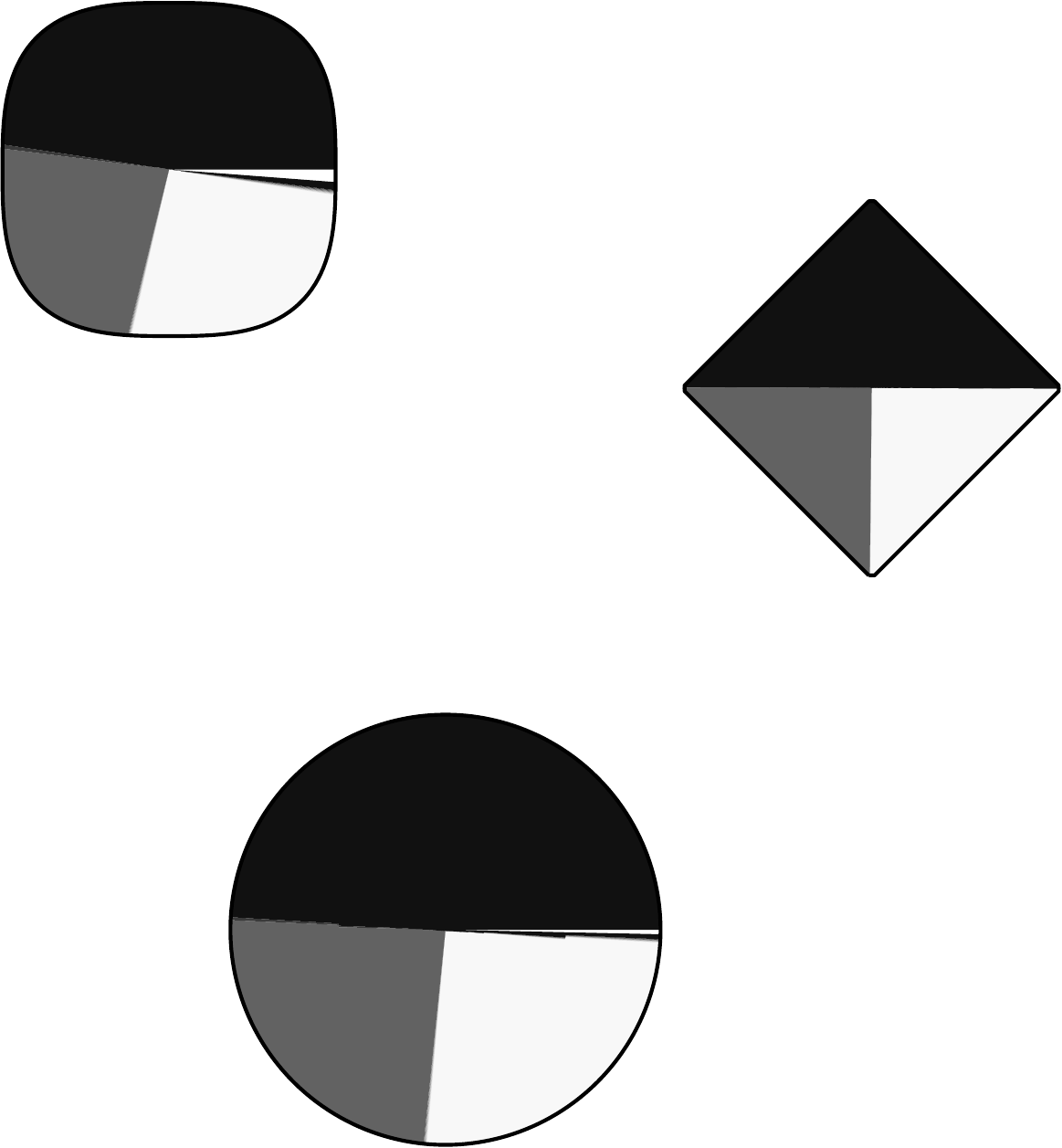}}
        \caption{CES preference function.}\label{fig:prefs_CES}
    \end{subfigure}~\begin{subfigure}[b]{.33\linewidth}
        \centering 
        \fbox{\includegraphics[width=0.8\linewidth]{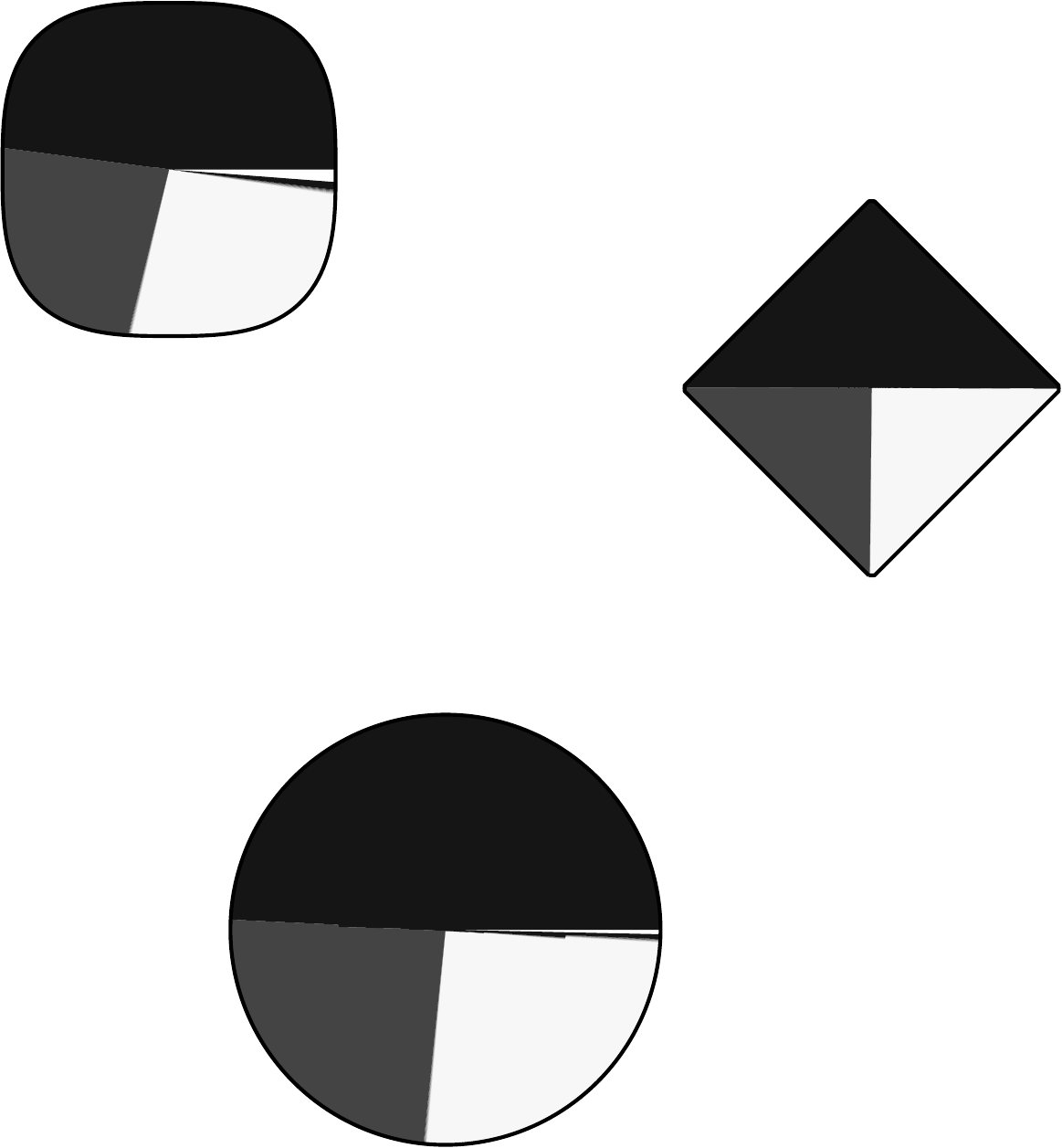}}
        \caption{LF preference function.}\label{fig:prefs_LF}
    \end{subfigure}
    \caption{Example of the three preference functions based on production models. From lightest to darkest, from lowest to highest minimum preference level.}\label{fig:prefs_prod}
\end{figure}

Given this methodology, Figure \ref{fig:prefs_prod} displays the three production-based preference functions for two factors ($m=2$). From left to right, we have the Cobb-Douglas (CD), the Constant Elasticity of Substitution (CES) utility function, and the Leontief function. Again, the preference level ranges from a lighter to a darker color, indicating the transition from the lowest to the highest minimum preference level for locating the entry point (from 0 to 1 in the normalized scale). 
In this example, the three functions appear quite similar because we used the same construction. However, this may not necessarily be the case in a real scenario. The only difference here is the preference level assigned to each subregion by each function.

\subsection{Thresholds}\label{sec:thresholds}

Note that the preference functions that we propose above are not required to take values in $[0,1]$, and the consideration of \textit{intuitive} thresholds in our models would require to normalize this function. The normalization that we propose is based on rescaling the preference functions as:
\begin{equation}
\widetilde{\Phi}_i(x) = \dfrac{\Phi_i(x) - LB(\Phi_i)}{UB(\Phi_i)-LB(\Phi_i)},\; \forall x \in S_i,
\end{equation}
\noindent where $UB(\Phi_i)= \max_{x \in S_i} \Phi(x)$ and $LB(\Phi_i) = \min_{x\in S_i} \Phi(x)$ (not that these values are achieved since the regions are assumed to be compact). 

The computation of $UB(\Phi_i)$ for all the (concave) functions that we propose requires solving a linear or a second order cone optimization problem (adequately rewritten), and then, solved in polynomial time. Calculating $LB(\Phi_i)$ is more challenging for the nonlinear preference functions since it requires solving a concave minimization problem over a convex region. We apply a spatial branch-and-bound procedure (that is actually implemented in some of the off-the-shelf software) to solve, exactly, these auxiliary problems.



\section{Numerical Experiments}\label{sec:experiments}

This section reports the results of our computational experience testing the different preference functions and justifying their use in location problems. 
For this purpose, we develop a completely random and independent dataset, which we describe in the Section \ref{sec:exp.data}.
Section \ref{sec:exp.perf} describes the computational performance of the problem with different functions in order to determine the computational complexity of the problem.
Finally, in Section \ref{sec:exp.quality} we make a detailed study of the solutions provided by the different preference functions.

All experiments were run on the Linux-based Supercomputer ``albaicin'', composed of 9520 Intel Xeon Cores at 2.7 GHz and 35 Tbytes RAM total, as 170 nodes interconnected with Infiniband HDR200 network. 
We use only 8 threads and 32 gigabyte of RAM for our experiments. The models were coded in Python 3.8.2 wich Gurobi 10.0.3 as optimization solver.
A time limit of 1 hour is fixed for all the instances.

\subsection{Instance generation} \label{sec:exp.data}

Due to the novelty of the approach proposed in this work, we have generated fully randomized instances to show the validity of using different preferences within regions. Thus, in this section, we detail how we generate the instances and the different parameters, so that the interested reader can replicate the experiments. All these instances are available in our GitHub repository:
\href{https://github.com/vblancoOR/preference_location}{\texttt{preference\_location}}.
The section is organized starting with the generation of the regions and then detailing the parameters. 

\begin{figure}[ht]
    \begin{subfigure}[b]{.5\linewidth}
        \centering 
        \fbox{\includegraphics[scale=0.3]{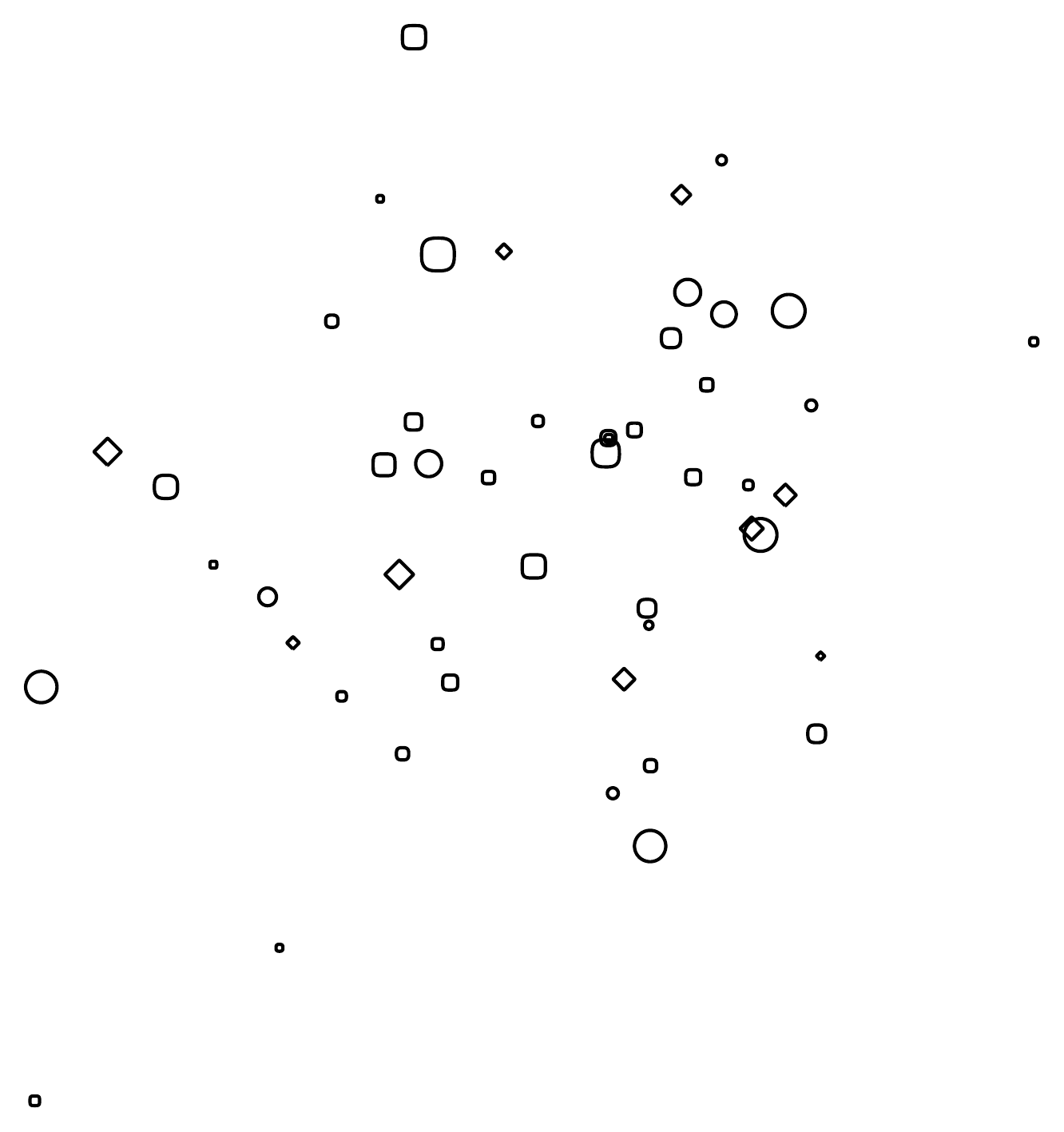}}
        \caption{$n=50$.}\label{fig:data50}
    \end{subfigure}~\begin{subfigure}[b]{.5\linewidth}
        \centering 
        \fbox{\includegraphics[scale=0.3]{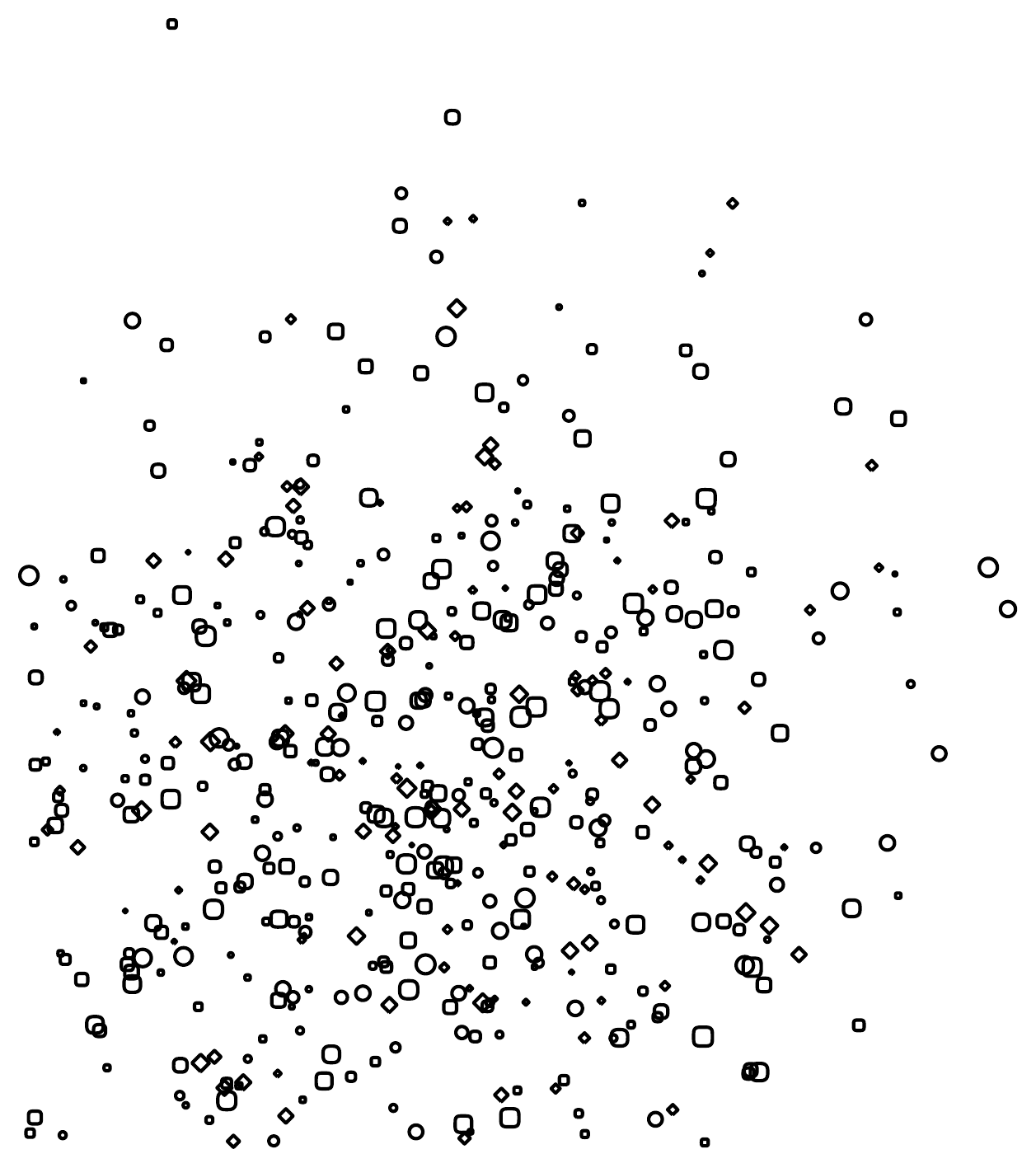}}
        \caption{$n=500$.}\label{fig:data500}
    \end{subfigure}~
    \caption{Dataset examples for same seed with $n=50$ and $n=500$.}\label{fig:data}
\end{figure}

We have generated five independent instances with $n$ regions ranging in $\{10,20,50,100,200,500\}$.   We generated planar $\ell_\tau$--norm balls regions with random radii and centers. Thus, each region is identified by its radius,  its center, its weight in the objective function, and the norm. The radii and the weights are chosen uniformly randomly in $(0.05,0.25)$ and $(0,1)$, respectively. We use \texttt{make\_blobs} function of the \texttt{sklearn} python library for generating the centers. Finally, we propose three different scenarios for both the norms and the distances between the regions and the facilities:
\begin{description}
\item[Scenario \texttt{l1}:] Both the distances and the balls are based on the $\ell_1$-norm. The objective function and the constraints assuring that the entry points are in the regions can be rewritten as linear contraints.
\item[Scenario \texttt{l2}:] 
Both the distances and the balls are based on the Euclidean norm. 
\item[Scenario \texttt{mixed}:] We randomly assign a norm $\ell_\tau$ with $\tau \in \{1,2,3,4\}$ to measure both the transportation cost and the ball defining the region. This scenario allow us to show that our proposal is valid for different distances, which may have applications in situations where the transportation modes are different for the regions~\citep[see e.g.]{blanco2017continuous}.
\end{description}
Figure \ref{fig:data} depicts an example of one of the sets considered for  sizes $50$ and $500$.

Three values of installed facilities are chosen with $p \in \{1,2,5\}$. By the computational difficulty of multi-facility problems, we solved the single-facilty problems for all sizes, and the multifacility instances for sizes $n\in \{10,20,50\}$.

Regarding the preferences, all functions shown in this paper are used for our experiments. We denote as \texttt{L}---Linear, \texttt{D}---Distance, \texttt{CES}---Constant Elasticity of Substitution, \texttt{CD}---Cobb-Douglas and \texttt{LF}---$m$-factor Leontief. 
Table \ref{tab:functions} shows the index and parameters involved in each function.

\begin{table}[ht]
    \begin{center}
        \begin{tabular}{ccl}
            \toprule
            Function & Index & \multicolumn{1}{c}{Parameters}\\
            \hline
            \texttt{L} & & $\gamma_i,\; \gamma_0,\; \forall i \in I$\\
            \texttt{D} & $K_i$, number of reference points & $\lambda_{ij},\; B_{ij},\; \forall i \in I, j \in \{1,\ldots,K_i\}$\\
            \texttt{CES} & $m$ number of factors & $\tau_i,\; \beta_{ij},\;  g_{ij},\; \forall i \in I, j \in \{1,\ldots,m\}$ \\
            \texttt{CD}, \texttt{LF}& $m$ number of factors & $\beta_{ij},\;  g_{ij},\; \forall i \in I, j \in \{1,\ldots,m\}$ \\
            \bottomrule
        \end{tabular}
        \caption{Different index and parameters involved in each function.}
        \label{tab:functions}
    \end{center}
\end{table}

The values that define the function \texttt{L} are chosen by the following way where $X$ is a random integer value in $[0,1]$ that defines the preference orientation within the region, and $Y,\tilde{Y}$ are two different random integer values in $[1,10]$:
$$
\gamma_0 = \tilde{Y},\; \gamma_i = (-1)^X Y, \forall i \in I.
$$

In case the function \texttt{D} is chosen, the number of reference points, $K_i$, is chosen as an integer random value in [1,3], for each region $i \in I$. The position of these points is randomly decided around the center defining the region. The $\lambda$--values are fixed for all regions $i \in I$, as $\lambda_{ij} = \frac{1}{K_i}, \forall j \in \{1,\ldots,K_i\}$.

\begin{figure}[ht]
    \begin{subfigure}[b]{.5\linewidth}
        \centering 
        \fbox{\includegraphics[scale=0.2]{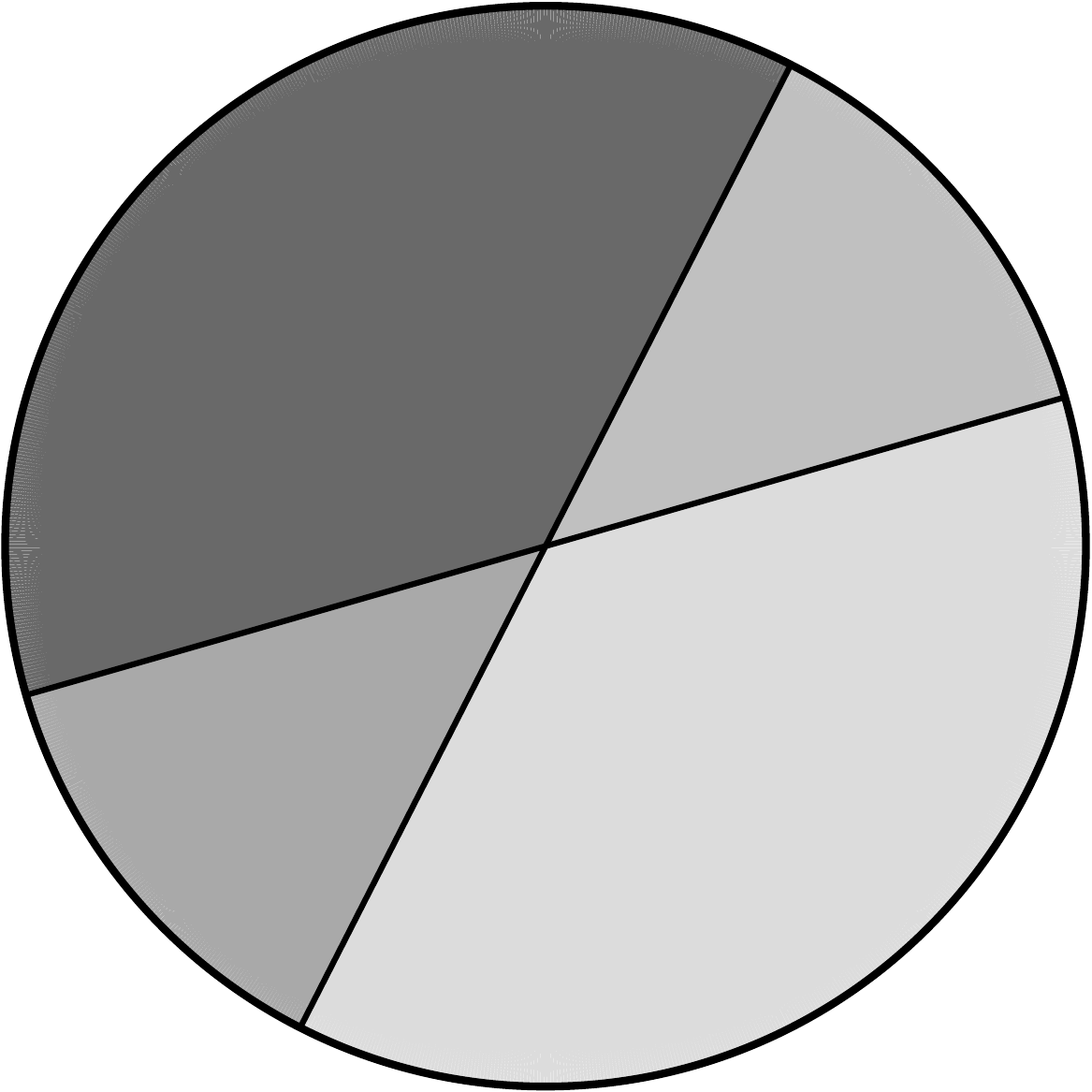}}
        \caption{First factor.}\label{fig:subregions1}
    \end{subfigure}~\begin{subfigure}[b]{.5\linewidth}
        \centering 
        \fbox{\includegraphics[scale=0.2]{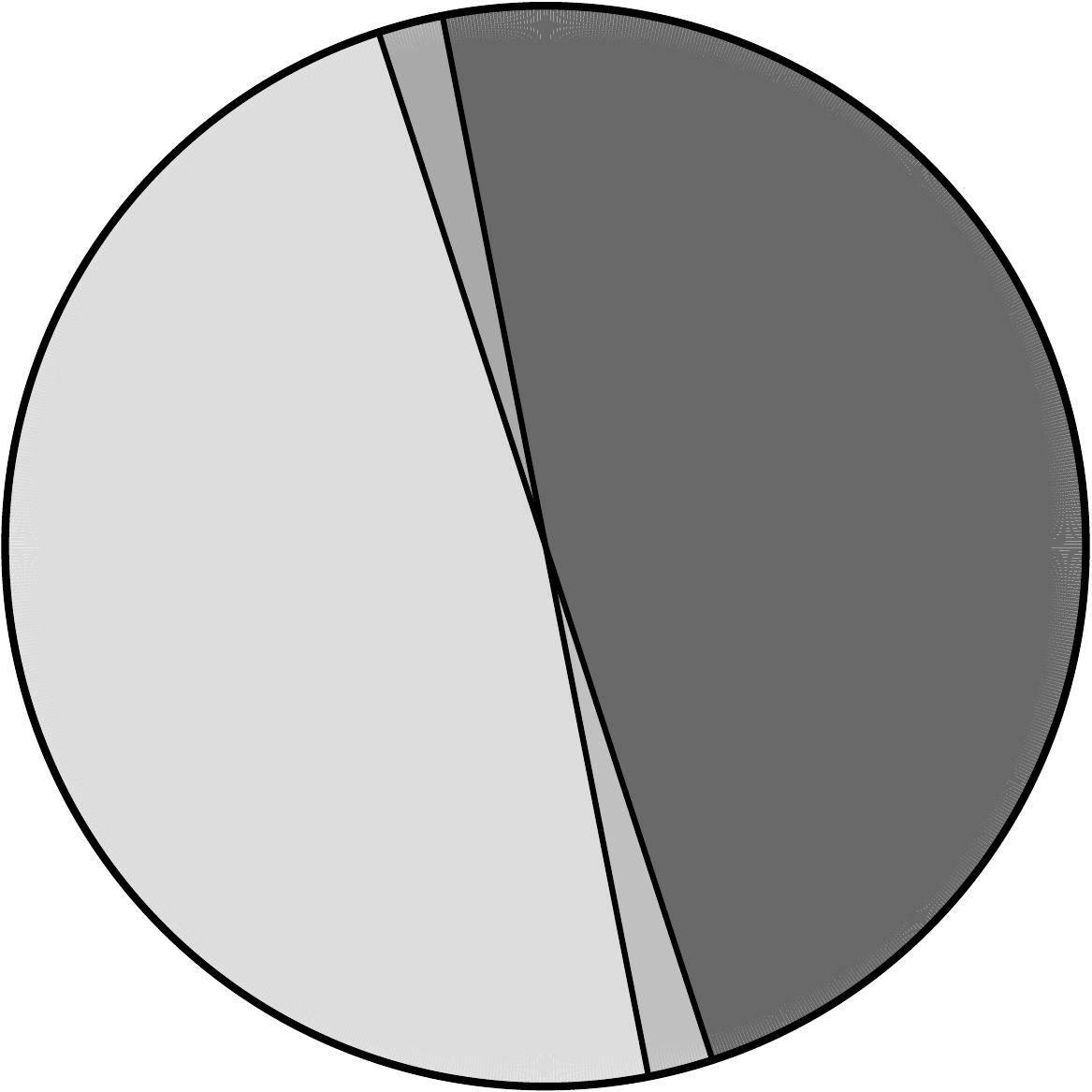}}
        \caption{Second factor.}\label{fig:subregions2}
    \end{subfigure}~
    \caption{Given two factors, $m=2$. Example of subdivision of the region into four subcells for $\ell_2$-norm ball.}\label{fig:subregions}
\end{figure}

For any preference based on production models, \texttt{CES}, \texttt{CD}, \texttt{LF}, we set the dimension of the production function, $m = 2$. 
This value indicates how many factors are considered in each region to define preference. Where for each factor $j$, the function $g$ will return the value of that preference. As already mentioned in Section \ref{sec:prefs}, we assume that for each factor $j$, the region is divided into cells by hyperplanes. The number of hyperplanes is set to $2$ resulting in a partition of the ball-shaped region into 4 sub-regions. These sub-regions are defined by two angles (randomly chosen in $[0,360]$). Figure \ref{fig:subregions} shows an example of how these regions look like. Finally, the $\beta$--value is fixed as $\beta = \{0.5,0.5\}$ and, in case \texttt{CES} function is selected, $\tau$ is fixed to $0.5$ for all $i \in I$. 

The preference functions were rescaled to $[0,1]$ as indicated in Section \ref{sec:thresholds}.  Three different thresholds are considered, $\Phi_i^* \in \{0,0.2,0.8\},\; \forall i \in I$. The case $\Phi_i^*=0$ coincides with the case in which the preferences are not considered and only the transportation costs has an impact in the solution. 

Finally, to show the performance of considering collocation (C-RLPP) or not (RLPP), we run the experiments with both options for the single facility case and we only run the experiment with collocation for the multi-facility case. In summary, 2700 instances were generated.

\subsection{Computational performance} \label{sec:exp.perf}

This section addresses the computational complexity of the models shown in this paper. We start with the case of one facility and then with the multi-facility case.

\subsubsection*{Single-facility}

The single facility case, as already mentioned, can be solved in polynomial time (the model can be cast as a continuous second order cone optimization problem for the non-collocation problem), so we are able to solve instances of larger sizes. Only $2$ out of the 2700 instances were not optimally solved within the time limit. These instances belong, one, to the scenario  without collocation, for the \texttt{CD} function  and $n=500$, and the other to the case with collocation, for the \texttt{CES} function and $n=500$. Given the very low percentage of unsolved instances, we excluded these two cases from the analysis in this section and will discuss them at the end.

\begin{figure}[ht]
    \begin{subfigure}[b]{.5\linewidth}
        \centering 
        \includegraphics[scale=0.35]{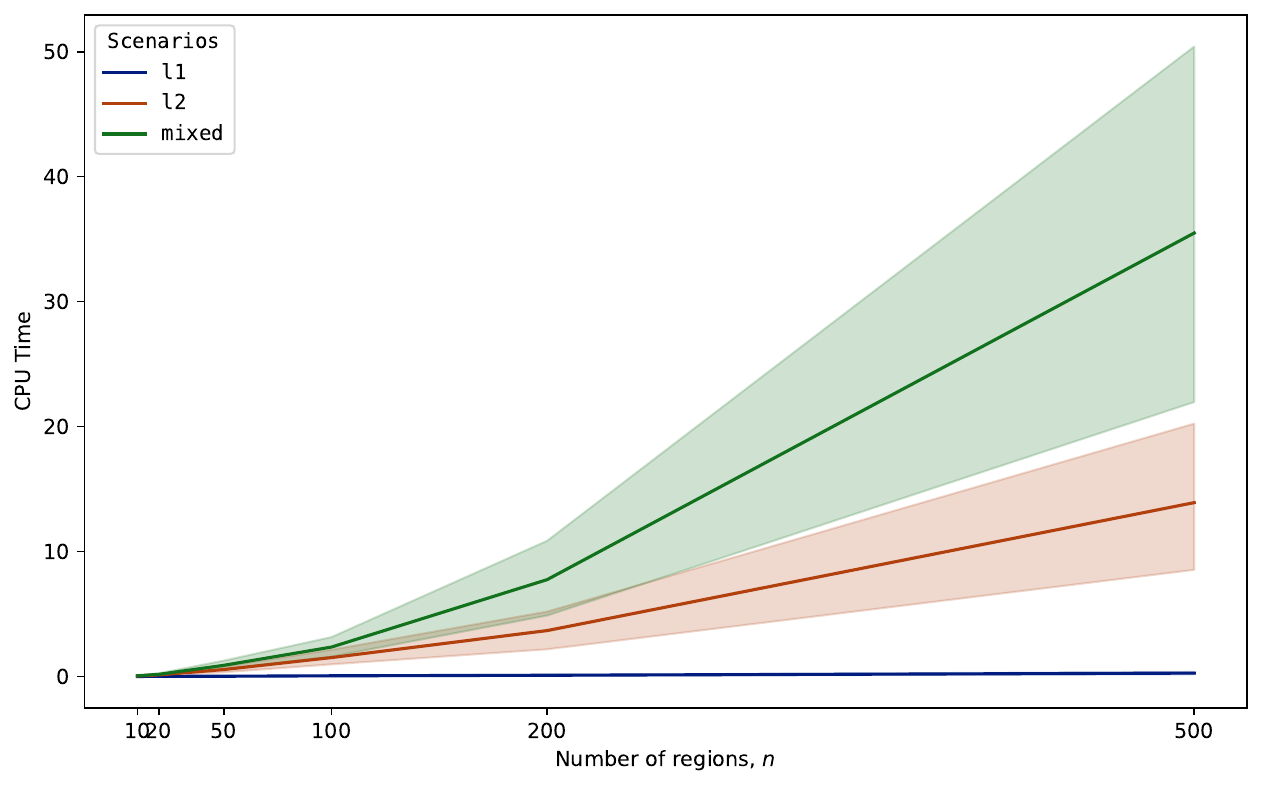}
        \caption{(RCLPP).}\label{fig:time_n_nrs0}
    \end{subfigure}~\begin{subfigure}[b]{.5\linewidth}
        \centering 
        \includegraphics[scale=0.35]{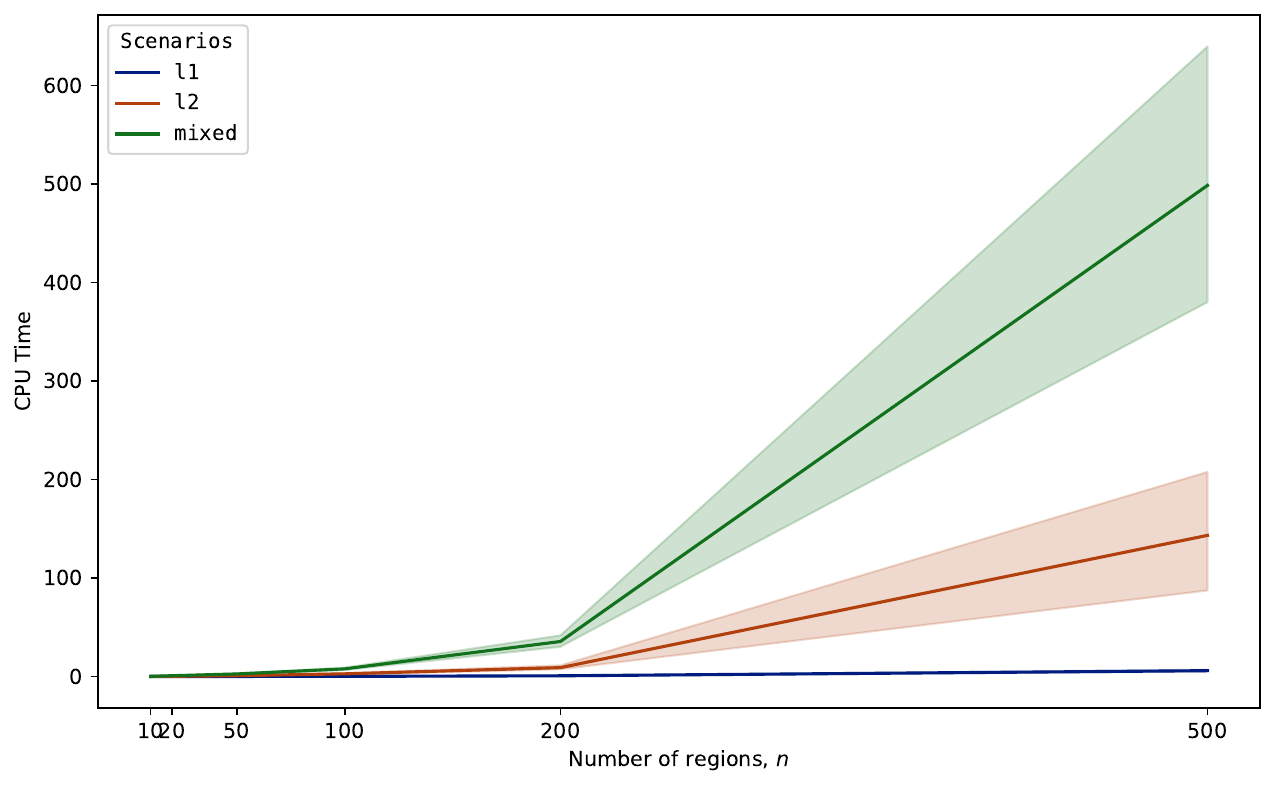}
        \caption{(C-RCLPP).}\label{fig:time_n_nrs1}
    \end{subfigure}~
    \caption{Average CPU times required for solving (RCLPP) and (C-RCLPP) by number of regions.}\label{fig:time_n_nrs}
\end{figure}

Figure \ref{fig:time_n_nrs} shows the average required CPU times as a function of dataset size for each of the scenarios in the study (norms), for both (RCLPP)--Figure \ref{fig:time_n_nrs0} and (C-RCLPP)--Figure \ref{fig:time_n_nrs1}. As expected, when using only the $\ell_1$ norm, the computational times required to solve the problems are significantly smaller than the other two scenarios. The mixed norm scenario contains the most challenging instances, being the required CPU time for solving the problem significantly larger than for the Euclidean norm scenario. This difference is even more evident when solving the problem with collocation (C-RCLPP) where additional binary variables are incorporated to the model.

\begin{figure}[ht]
 \begin{subfigure}[b]{.5\linewidth}
\centering 
\includegraphics[scale=0.35]{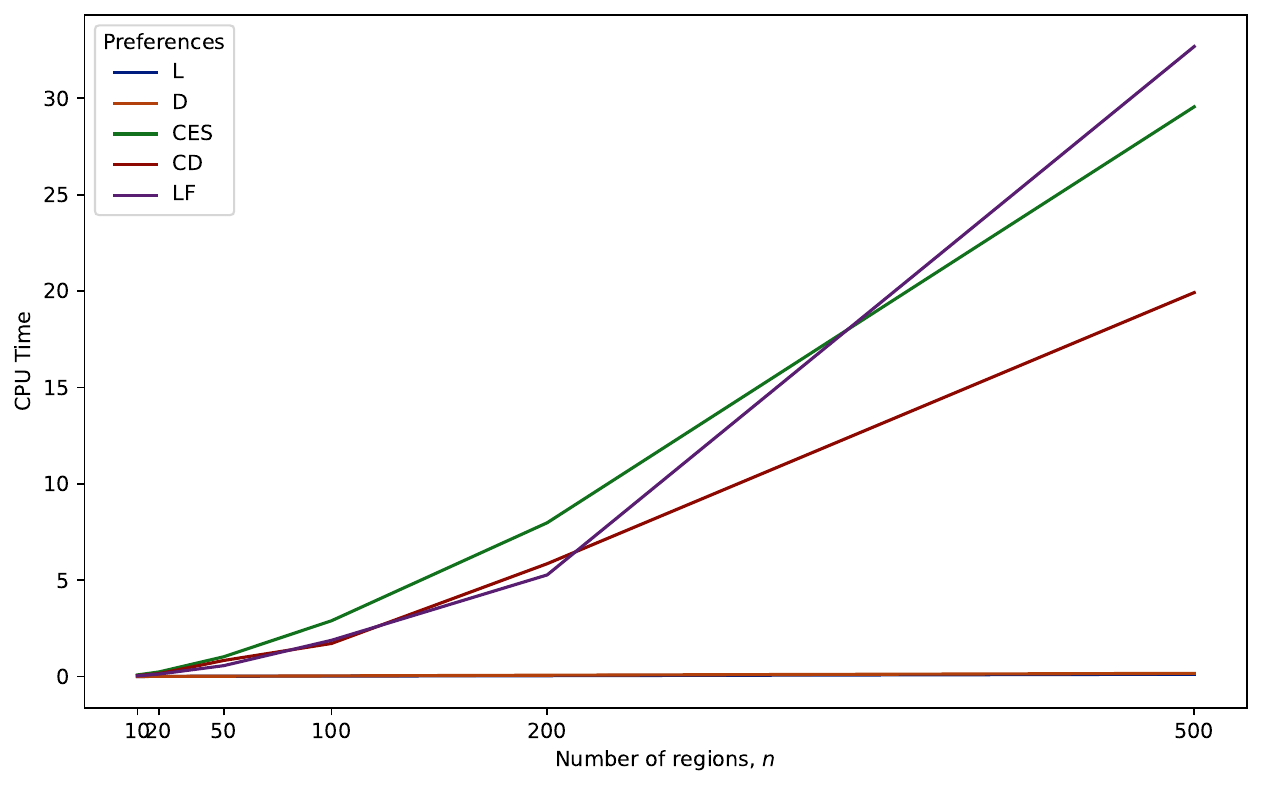}
\caption{(RCLPP).}\label{fig:time_n_ff0}
 \end{subfigure}~\begin{subfigure}[b]{.5\linewidth}
\centering 
\includegraphics[scale=0.35]{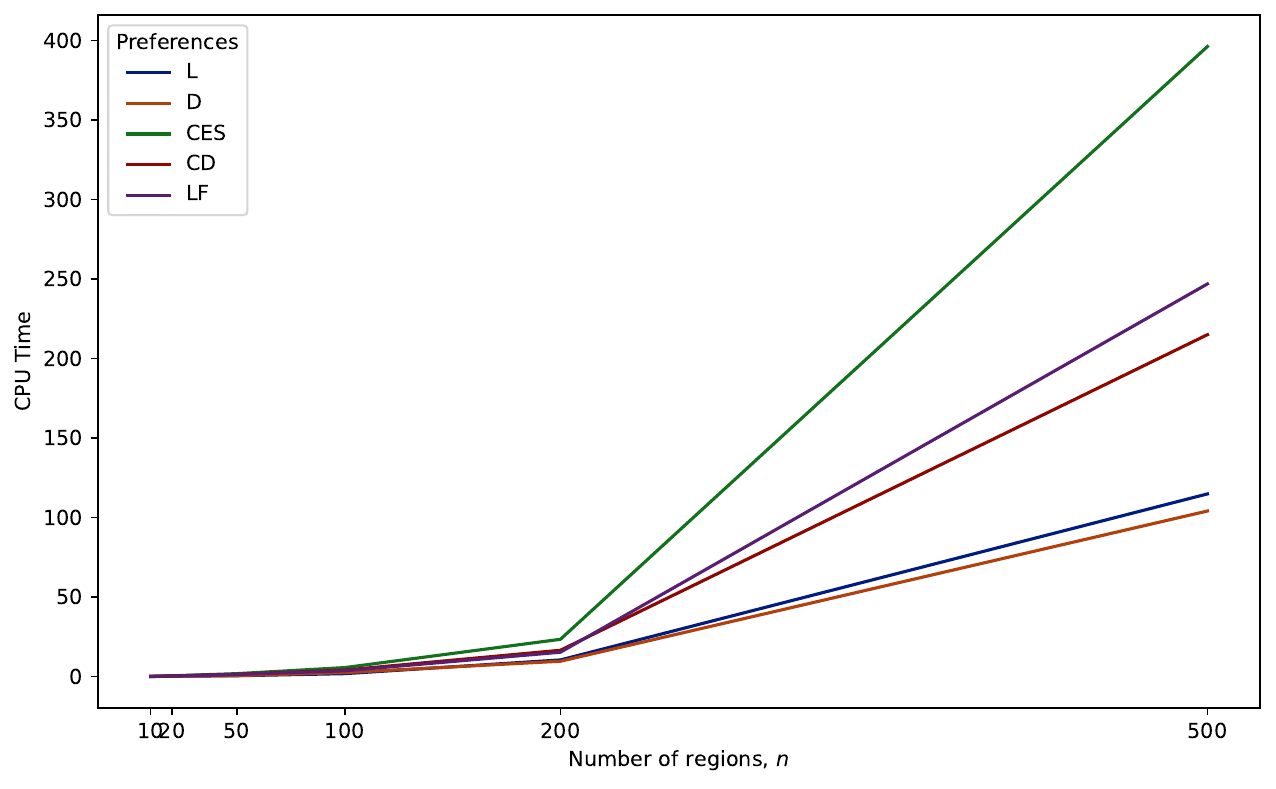}
\caption{(C-RCLPP).}\label{fig:time_n_ff1}
 \end{subfigure}~
 \caption{Average computational times as a function of dataset size, by the preferences functions used in the study and differentiating on considering collocation or not.}\label{fig:time_n_ff}
\end{figure}

The differences, in CPU time, by the preference functions used in the models are reported in Figure \ref{fig:time_n_ff}. For (RCLPP),  the linear and distance based preference functions only required a few seconds to be optimally solved, whereas the preferences based on production models required more, but still reasonable, CPU time. For (C-RCLPP), CPU times dramatically increase by the number of binary variables. While all three production-based preference functions have a similar performance for (RCLPP), it seems that \texttt{CES} is more challenging for (C-RCLPP).

\begin{figure}[ht]
 \begin{subfigure}[b]{.5\linewidth}
\centering 
\includegraphics[scale=0.35]{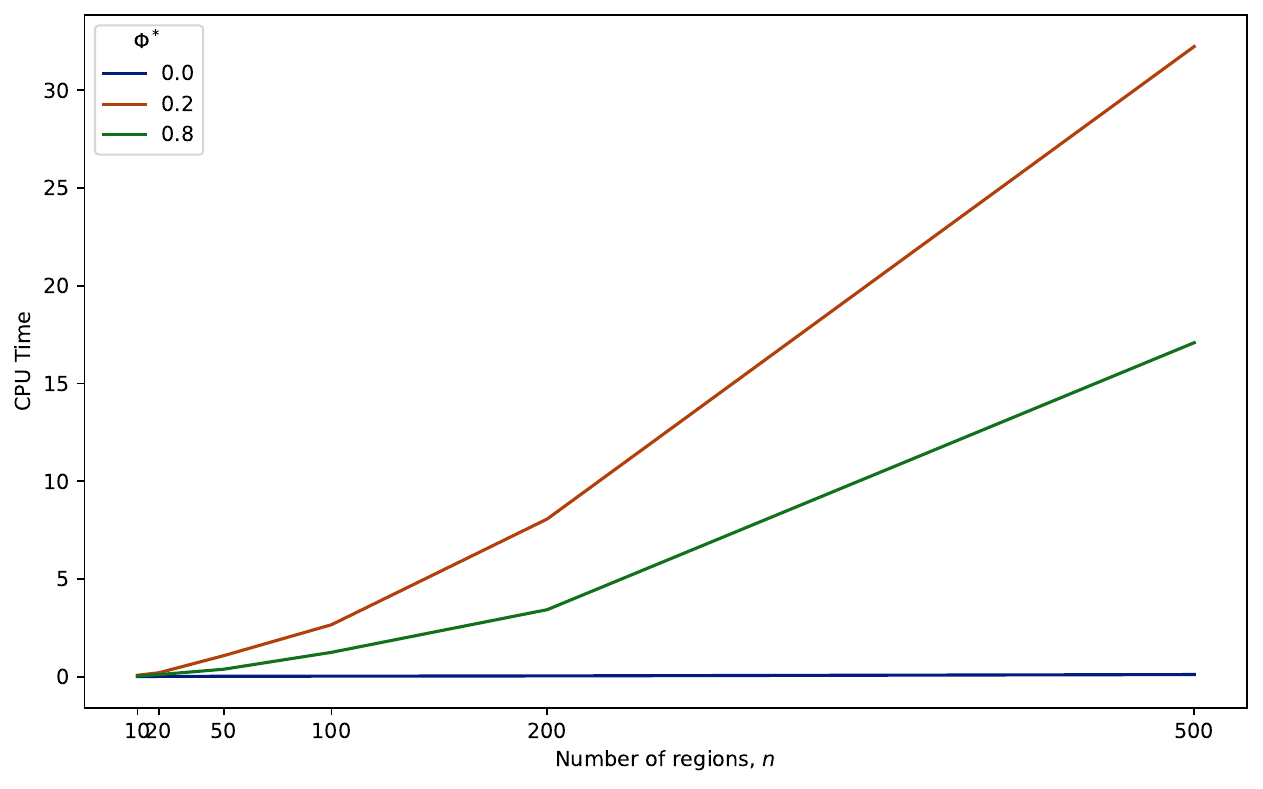}
\caption{(RCLPP).}\label{fig:time_n_thr0}
 \end{subfigure}~\begin{subfigure}[b]{.5\linewidth}
\centering 
\includegraphics[scale=0.35]{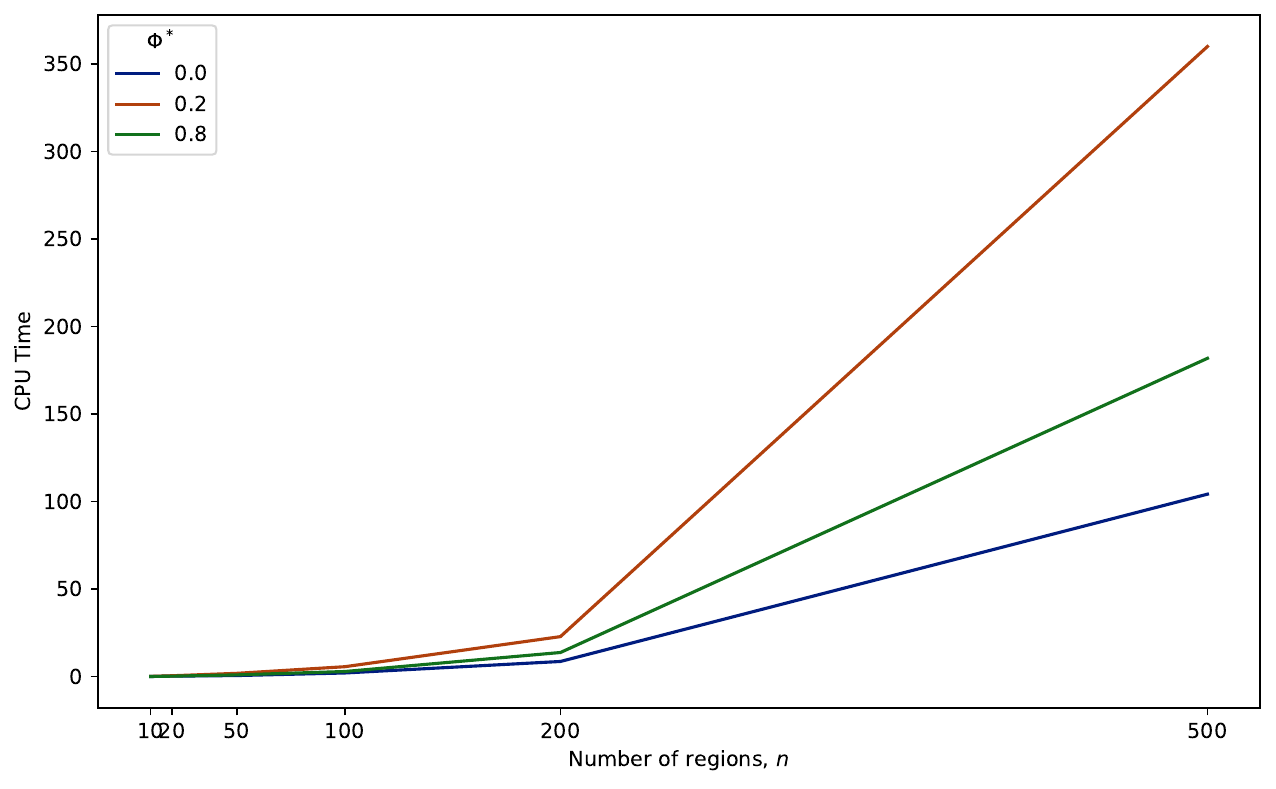}
\caption{(C-RCLPP).}\label{fig:time_n_thr1}
 \end{subfigure}~
 \caption{Average computational times as a function of dataset size, divided by each threshold used in the study and differentiating on considering collocation or not.}\label{fig:time_n_thr}
\end{figure}

Figure \ref{fig:time_n_thr} shows the average CPU times based on size, by the thresholds used. For $\Phi^*=0$, the problem simplifies to a neighborhood version of the Weber Problem, and it can be modeled as a continuous second order cone optimization problem that can be solved in polynomial time. The most computationally challenging for the threshold $0.2$, which require more CPU time than the case $0.8$. This is due to the reduction of the regions $P_i$. While the case $0.2$ requires searching in the whole region (perhaps except the optimal solution obtained when the threshold is $0$) for the optimal position of the entry points, the $0.8$ case is a more restrictive area closer to the ideal preference for the users.

Finally, in Table \ref{tab:time_500} we further analyze the most challenging size, $n=500$. In this table, we show the minimum, average, and maximum CPU times required to optimally solve the instances, classified by dataset size ($n$), norm used (scenario), and threshold ($\Phi^*$). The thresholds are only disaggregated for scenario \texttt{mixed}, where we observed differences between the different values. Each row with a specified threshold includes data from 5 instances (one per seed), while the rest cover 15 instances.

One can observe from the table that for any preference function, the model can solve 500 demand sets for scenarios \texttt{L1} and \texttt{L2}. However, the \texttt{mixed} scenario is more computationally challenging. Note that when rewriting as SOC constraints the expressions (objective function and constraint) where the $\ell_3$ and $\ell_4$ appear, require incorporating more SOC constraints than for the Euclidean norm, and then increasing its difficulty, specially for thresholds $\Phi^* > 0$. As already observed, our results indicate that the instances for $\Phi^* = 0.2$ require more CPU time to be optimally solved than those where $\Phi^* = 0.8$.

\begin{table}[h!]
\centering
\adjustbox{width=\textwidth}{
\begin{tabular}{@{}lllrrrrrrrr@{}}
\toprule
Preferences& \texttt{scenario} & $\Phi^*$ & \multicolumn{4}{c}{(RCLPP)}& \multicolumn{4}{c}{(C-RCLPP)}\\ 
\cline{4-7}\cline{8-11}
&&  & \multicolumn{1}{c}{Min} & \multicolumn{1}{c}{Average} & \multicolumn{1}{c}{Max} & \multicolumn{1}{c}{\texttt{TL}} & \multicolumn{1}{c}{Min} & \multicolumn{1}{c}{Average} & \multicolumn{1}{c}{Max} & \multicolumn{1}{c}{\texttt{TL}} \\
\hline
\texttt{L}& \texttt{L1} &  & 0.03 & 0.04& 0.05 && 3.56 & 4.24& 5.89 &\\
\cline{4-11}
& \texttt{L2} &  & 0.09 & 0.10 & 0.12 && 13.12& 55.62 & 132.06  &\\
\cline{4-11}
& \texttt{mixed} & 0& 0.18 & 0.19 & 0.20 && 189.52  & 229.45& 318.59  &\\
&& 0.2& 0.18 & 0.21& 0.23 && 169.16  & 448.99& 1307.92 &\\
&& 0.8& 0.21 & 0.28& 0.51 && 102.76  & 175.90& 384.89  &\\
\hline
\texttt{D}& \texttt{L1} &  & 0.05 & 0.06& 0.08 && 3.31 & 4.00& 5.83 &\\
\cline{4-11}
& \texttt{L2} &  & 0.08 & 0.14& 0.18 && 8.79 & 55.46 & 133.86  &\\
\cline{4-11}
& \texttt{mixed} & 0& 0.18 & 0.19& 0.20 && 188.34  & 229.75& 316.93  &\\
&& 0.2& 0.33 & 0.35& 0.39 && 219.49  & 389.89& 611.87  &\\
&& 0.8& 0.36 & 0.40& 0.45 && 126.31  & 139.09& 156.32  &\\
\hline
\texttt{CES} & \texttt{L1} &  & 0.04 & 0.63& 1.93 && 3.42 & 9.19& 18.85&\\
\cline{4-11}
& \texttt{L2} &  & 0.10 & 29.75 & 106.49  && 15.74& 292.05& 1564.61 &\\
\cline{4-11}
& \texttt{mixed} & 0& 0.18 & 0.19& 0.20 && 189.62  & 229.93& 319.73  &\\
&& 0.2& 25.44& 92.48 & 221.25  && 618.06  & 1418.05  & 2041.17 & 1 (0.01\%)\\
&& 0.8& 0.76 & 82.21 & 174.11  && 214.31  & 1218.98  & 3462.80 &\\
\hline
\texttt{CD}  & \texttt{L1} &  & 0.04 & 0.34& 1.28 && 3.31 & 7.14& 16.03&\\
\cline{4-11}
& \texttt{L2} &  & 0.09 & 14.44 & 65.35&& 9.91 & 92.98 & 213.59  &\\
\cline{4-11}
& \texttt{mixed} & 0& 0.18 & 0.19& 0.20 && 189.84  & 229.66& 319.07  &\\
&& 0.2& 51.89& 101.20& 213.79  && 274.16  & 1045.35  & 2918.93 &\\
&& 0.8& 0.59 & 36.97 & 82.05& 1 (99.98\%)  & 73.59& 359.97& 1014.26 &\\
\hline
\texttt{LF}  & \texttt{L1} &  & 0.04 & 0.30& 0.63 && 2.97 & 5.04& 8.14 &\\
\cline{4-11}
& \texttt{L2} &  & 0.10 & 25.16 & 99.31&& 13.24& 220.21& 1164.93 &\\
\cline{4-11}
& \texttt{mixed} & 0& 0.18 & 0.19& 0.21 && 190.90  & 276.65& 427.15  &\\
&& 0.2& 33.10& 138.42& 251.27  && 193.46  & 755.22& 1535.84 &\\
&& 0.8& 0.47 & 79.17 & 212.11  && 71.46& 515.05& 1529.31 &\\ 
\bottomrule
\end{tabular}
}
\caption{Results for the single-facility problem for $n=500$.} \label{tab:time_500}
\end{table}

To conclude, we discuss the case of the two instances that were not optimally solved within the time limit. These two instances are marked in columns \texttt{TL} where the number of unsolved instances in the row and the average MIPGap (in parenthesis) are indicated. Surprisingly, for (RCLPP) and  \texttt{CD} preference function with threshold $\Phi^* = 0.2$, the model failed to find any feasible solution that meets the threshold across all 500 regions within the time limit. After carefully analyzing the instance, we observed that this issue was caused by numerical precision within regions divided into cells; some cells have a minimum preference of 0.800001. As a result, the model searches within its tolerance limits for an entry point but fails to find a solution. In the other case, (C-RCLPP) and \texttt{CES} preference function, our approach was not able to certify optimality within the time limit, although its MIPGap is considerably small ($0.01\%$).

\subsubsection*{Multi-facility}

In the multi-facility case, we only report the results for $n \in \{10,20,50\}$ and model (C-RCLPP), for the sake of simplicity and analyze the limitations of our model. In this section, we present the same table but disaggregated for the different parameters: number of regions ($n$)-- Table \ref{tab:time_mf_n}, scenarios-- Table \ref{tab:time_mf_nrs}, thresholds ($\Phi^*$)-- Table \ref{tab:time_mf_threshold}, and preference function ($\Phi$)--Table \ref{tab:time_mf_function}. In all these tables we report the minimum, average, and maximum CPU Time for instances solved within the time limit. Additionally, each table includes two columns for instances that could not be solved within the time limit: the average MIPGap, the total number of unsolved instances for that combination, and the total number of instances per row in parenthesis.

From Table \ref{tab:time_mf_n} one can observe, as expected, that the model becomes more challenging as the number of regions increases. The model was able to solve all instances of size $n \in \{10,20\}$ when $p=2$, and for $p=5$ all the instances are solved when $n=10$. For $n=50$, the model was unable to solve $170$ out of $225$ instances within the time limit, meaning only about $25\%$ are solved. For $p=5$, none of the instances were optimally solved within the time limit, and the MIPGap after this time was, in average, $69.92\%$. 

Regarding computational times, we observe significant differences in the cases with unsolved instances, indicating that, although some instances were solved, they took a considerable amount of time.

\begin{table}[ht]
\centering
\adjustbox{width=\textwidth}{
\begin{tabular}{@{}llrrrrr@{}}
\toprule
$p$ & $n$ & Min CPU Time & Av. CPU Time & Max CPU Time & Av. MIPGap & \#TL (225) \\
\hline
2 & 10  & 0.04 & 0.49 & 5.50&  & 0 \\
 & 20  & 0.16 & 3.80 & 24.49  &  & 0 \\
 & 50  & 91.83 & 1016.28 & 3435.80 & 18.93 & 170  \\
 \hline
5 & 10  & 0.12 & 14.43  & 60.78  &  & 0 \\
 & 20  & 2.09 & 650.29 & 3595.60 & 28.12 & 80\\
 & 50  &  &  &  & 69.92 & 225 \\
 \bottomrule
\end{tabular}
}
\caption{Results for the multi-facility problems by number of facilities ($p$) and regions ($n$).} \label{tab:time_mf_n}
\end{table}

From Table \ref{tab:time_mf_nrs} we observe that unsolved instances appear in all the scenarios. The \texttt{mixed} scenario is significantly more difficult than the other scenarios. We observe also a high dispersion in times for the case $p = 2$ and scenario \texttt{L2} might seem striking if they are compared with \texttt{nrs} = 3. This is because a few of the instances were solved for $n=50$ for scenario \texttt{L2}, while none were solved when scenario \texttt{mixed}. In fact, the maximum coincides with that of Table \ref{tab:time_mf_n}.

\begin{table}[ht]
    \centering
    \adjustbox{width=\textwidth}{
    \begin{tabular}{@{}llrrrrr@{}}
        \toprule
        $p$ & \texttt{scenario} & Min CPU Time & Av. CPU Time & Max CPU Time & Av. MIPGap & \#TL (225) \\
        \hline
        2& \texttt{L1}& 0.04 & 248.75 & 2552.0 1& 8.23 & 22\\
         & \texttt{L2}& 0.10 & 37.47  & 3435.80 & 18.31 & 73\\
         & \texttt{mixed}& 0.25 & 4.46 & 24.49  & 22.68 & 75\\
         \hline
        5& \texttt{L1}& 0.12 & 6.06 & 45.99  & 41.65 & 75\\
         & \texttt{L2}& 6.38 & 440.84 & 3595.60 & 63.85& 110  \\
         & \texttt{mixed}& 14.29  & 437.48 & 3546.94 & 65.29& 120 \\
         \bottomrule
    \end{tabular}
    }
    \caption{Results for the multi-facility problems by number of facilities ($p$) and scenarios.} \label{tab:time_mf_nrs}
\end{table}

When dissagregating by threshold values, we observe from Table \ref{tab:time_mf_threshold} that the unsolved instances are not primarily caused by the incorporation of preferences $(\Phi^*=0$ versus $\Phi^*>0$). However, what can be seen in this table is that the average time to solve the instances increases slightly with the incorporation of preferences, and the maximum time also increases. (Note that the $p=5$ case is incomplete in terms of instances, as it does not solve any for $n=50$, which skews the maximum values.)

\begin{table}[ht]
    \centering
    \adjustbox{width=\textwidth}{
    \begin{tabular}{@{}llrrrrr@{}}
        \toprule
        $p$ & $\Phi^*$ & Min CPU Time & Av. CPU Time & Max CPU Time & Av. MIPGap & \#TL (225) \\
        \hline
        2& 0 & 0.05 & 94.83  & 1564.11 & 17.31 & 55\\
         & 0.2  & 0.05 & 105.78 & 2504.66 & 19.75 & 59\\
         & 0.8  & 0.04 & 137.16 & 3435.80 & 19.66 & 56\\
        \hline
        5& 0 & 0.12 & 219.13 & 3595.60 & 57.61 & 101  \\
         & 0.2  & 0.14 & 248.31 & 2706.46 & 59.39 & 104  \\
         & 0.8  & 0.12 & 322.58 & 3546.94 & 59.86 & 100 \\
         \bottomrule
\end{tabular}
}
\caption{Results for the multi-facility problems by number of facilities ($p$) and thresholds ($\Phi^*$).} \label{tab:time_mf_threshold}
\end{table}

Finally, in Table \ref{tab:time_mf_function}, when disaggregating by preference functions, the number of unsolved instances is roughly similar for each function, showing that in the multi-facility case, which already involves solving a mixed-integer problem, the complexity of solving remains about the same regardless of the preference function used. Additionally, we observe that the time dispersion behaves similarly across functions, and the average times are also comparable, except for the \texttt{CD} case with $p=5$.

\begin{table}[ht]
    \centering
    \adjustbox{width=\textwidth}{
    \begin{tabular}{@{}llrrrrr@{}}
        \toprule
        $p$ & $\Phi$ & Min CPU Time & Av. CPU Time & Max CPU Time & Av. MIPGap & \#TL (135) \\
        \hline
        2& \texttt{L}  & 0.05 & 97.75  & 1659.39 & 17.86 & 33\\
         & \texttt{D}  & 0.05 & 144.58 & 3435.80 & 19.13 & 32\\
         & \texttt{CES} & 0.05 & 99.51  & 2504.66 & 18.90 & 36\\
         & \texttt{CD} & 0.04 & 102.44 & 2064.14 & 19.07 & 35\\
         & \texttt{LF} & 0.04 & 117.86 & 2552.01 & 19.68 & 34\\
         \hline
        5& \texttt{L}  & 0.14 & 195.56 & 3592.30 & 57.42 & 62\\
         & \texttt{D}  & 0.12 & 273.80 & 3595.60 & 58.80 & 62\\
         & \texttt{CES} & 0.12 & 330.84 & 3546.94 & 60.76 & 59\\
         & \texttt{CD} & 0.12 & 186.03 & 2343.51 & 57.55 & 63\\
         & \texttt{LF} & 0.12 & 325.52 & 3422.14 & 60.45 & 59  \\
         \bottomrule
    \end{tabular}
    }
    \caption{Results for the multi-facility problems by number of facilities ($p$) and preference functions ($\Phi$).} \label{tab:time_mf_function}
\end{table}

In conclusion, regarding the computational times for the multi-facility case, times increase drastically for smaller sizes, with sets like $n=50$ and $p=5$ remaining unsolved. Using any preference function and the minimum preference in the regions does not significantly affect computational times. However, the type of norm used does influence time, with those involving fewer SOC constraints being easier to solve.

\subsection{Managerial insights: The impact of customers' satisfaction into the efficiency} \label{sec:exp.quality}

In this section, we analyze the question:  Which is the impact of the incorporation of preferences in the efficiency of the solution? To answer this question we conduct an exploratory analysis of the solutions obtained with our approaches. We compare the objective values obtained with the different thresholds for the preference functions (in particular when the threshold is zero is equivalent to the case when no preferences are considered, and then, the \textit{most efficient} solution.

In this analysis, we focus on the single-facility case, where a larger sample of information was obtained in our experiments: $2700$ instances solved within the time limit (excluding 2 instances that were removed from this analysis). The study of the multifacility case will require the development of solution methods capable of addressing a larger number of problems and will be the focus of future research in this area.

\begin{figure}[ht]
    \begin{subfigure}[b]{.5\linewidth}
        \centering 
        \includegraphics[scale=0.32]{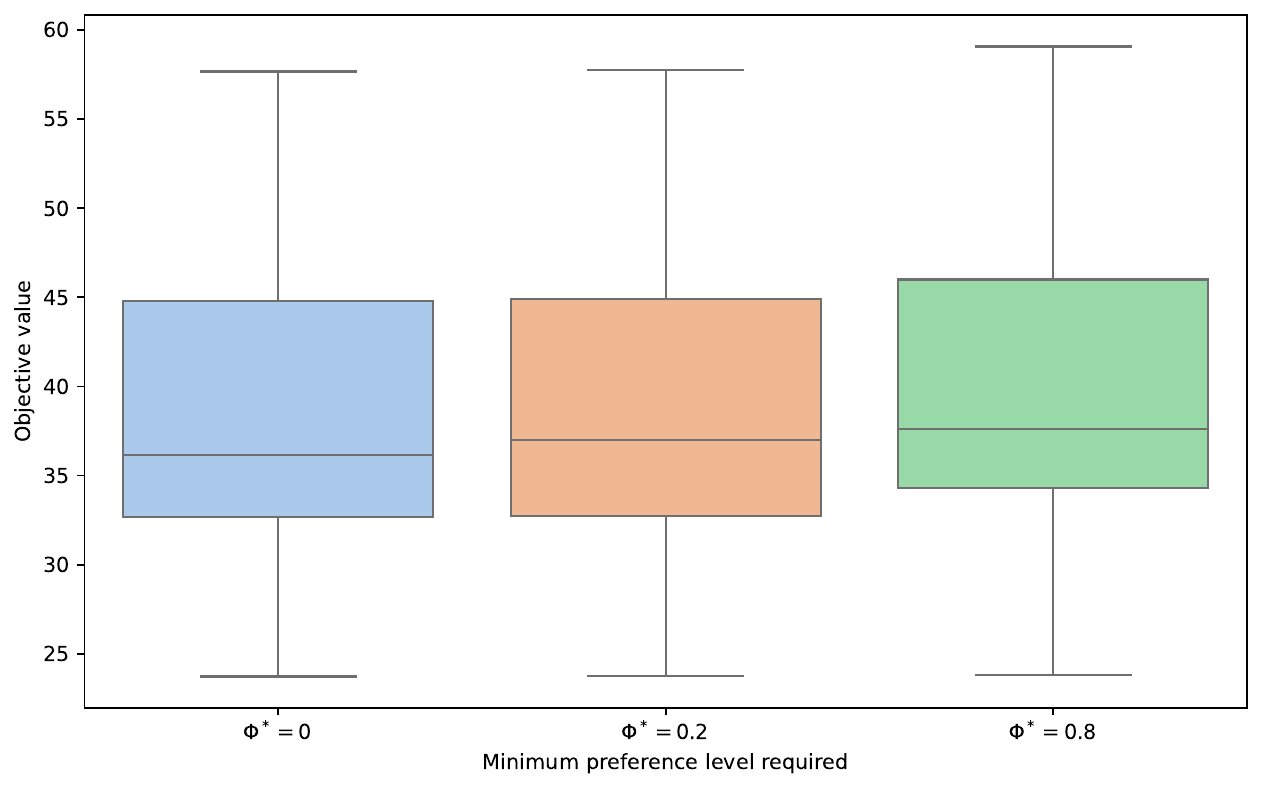}
        \caption{$n=20$.}\label{fig:boxplot20}
    \end{subfigure}~\begin{subfigure}[b]{.5\linewidth}
        \centering 
        \includegraphics[scale=0.32]{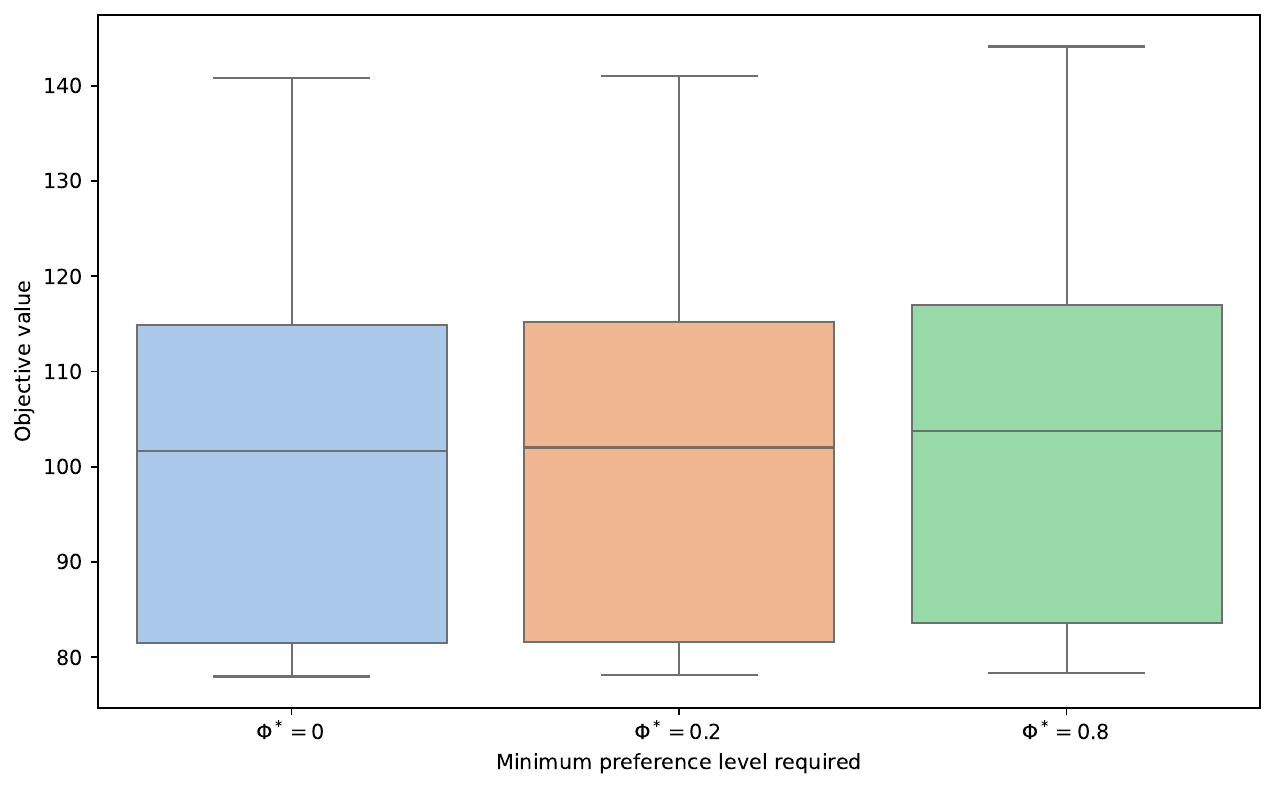}
        \caption{$n=50$.}\label{fig:boxplot50}
    \end{subfigure}\\ \begin{subfigure}[b]{.5\linewidth}
        \centering 
        \includegraphics[scale=0.32]{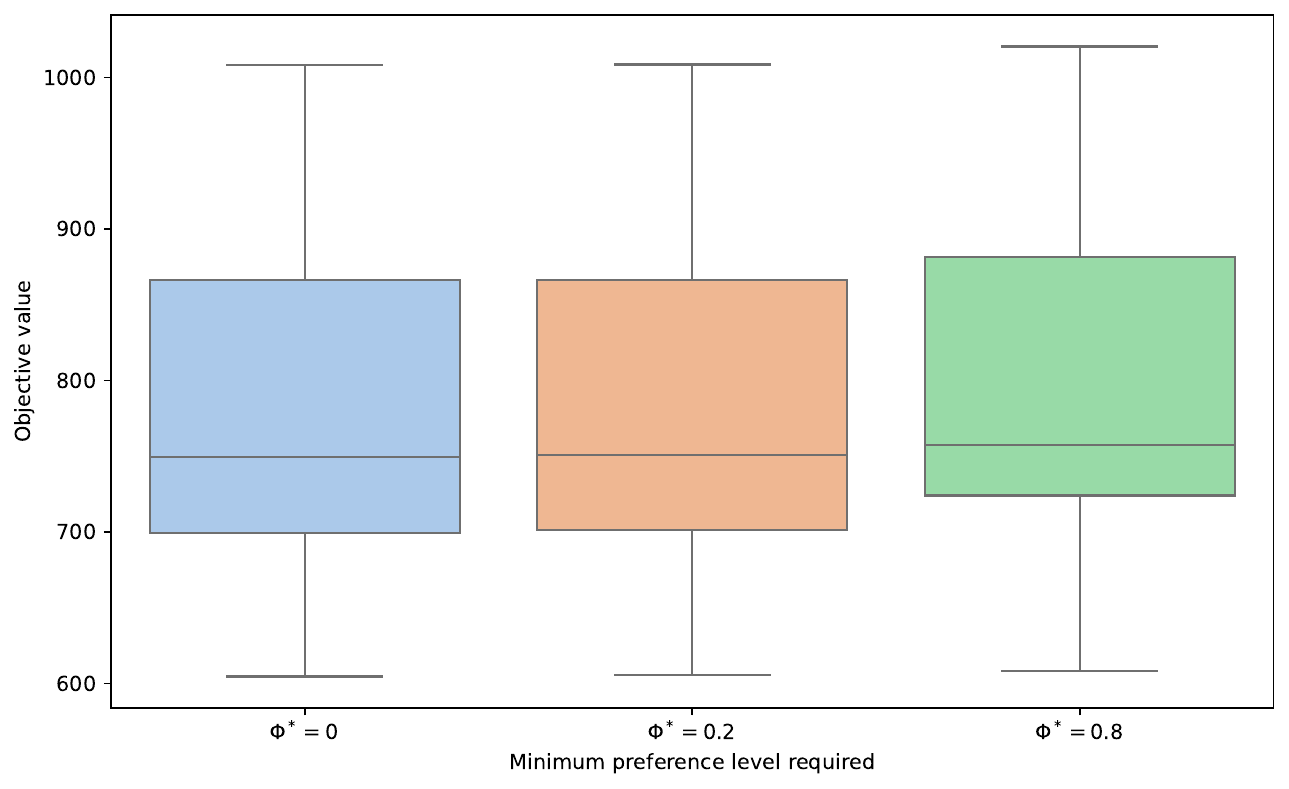}
        \caption{$n=200$.}\label{fig:boxplot200}
    \end{subfigure}~\begin{subfigure}[b]{.5\linewidth}
        \centering 
        \includegraphics[scale=0.32]{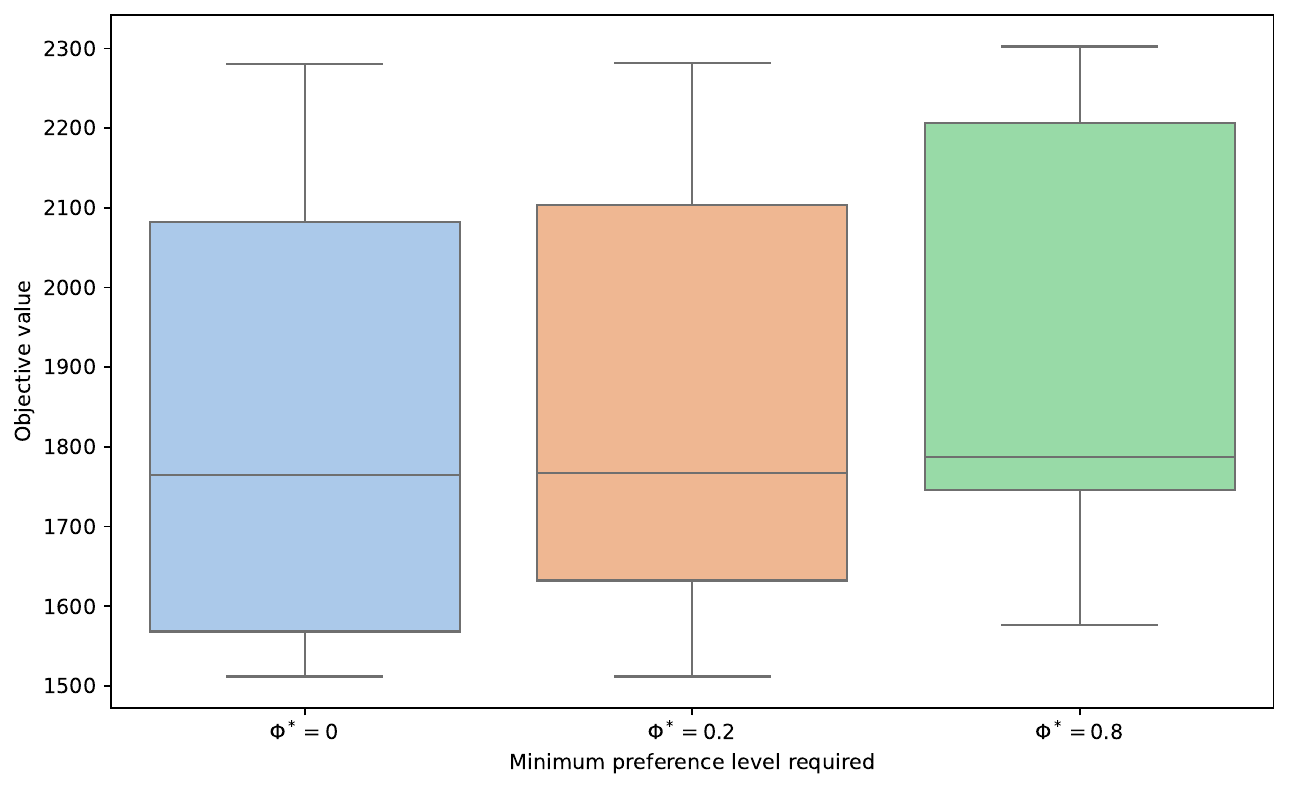}
        \caption{$n=500$.}\label{fig:boxplot500}
    \end{subfigure}~
    \caption{Boxplots of the objective value comparing with the preference level ($\Phi^*$) for different values of $n$.}\label{fig:boxplot}
\end{figure}

To facilitate the exposition we draw several boxplots in Figure \ref{fig:boxplot} to compare the objective values (transportation costs) obtained for different number of regions, $n \in \{20,50,200,500\}$, by preference thresholds ($\Phi^*$). From these boxplots, we observe that the median of transportation costs in the optimal solutions does not change notably when comparing the three groups, being the difference more distinguishable in the case of $n=500$. 
The differences become more pronounced when analyzing the shape of the boxplots. Specifically, for \( n=20 \), notable changes in the interquartile range (IQR) distributions are observed when comparing \(\Phi^* = 0\) to \(\Phi^* = 0.8\). In contrast, for \( n=500 \), these changes are apparent across all values of \(\Phi^*\).
Regardless, small changes in the objective value are advantageous for our approach, as it indicates that the additional costs for the decision-maker when considering customer preferences is not substantially higher compared to ignoring them.

\begin{figure}[ht]
    \begin{subfigure}[b]{.5\linewidth}
        \centering 
        \includegraphics[scale=0.32]{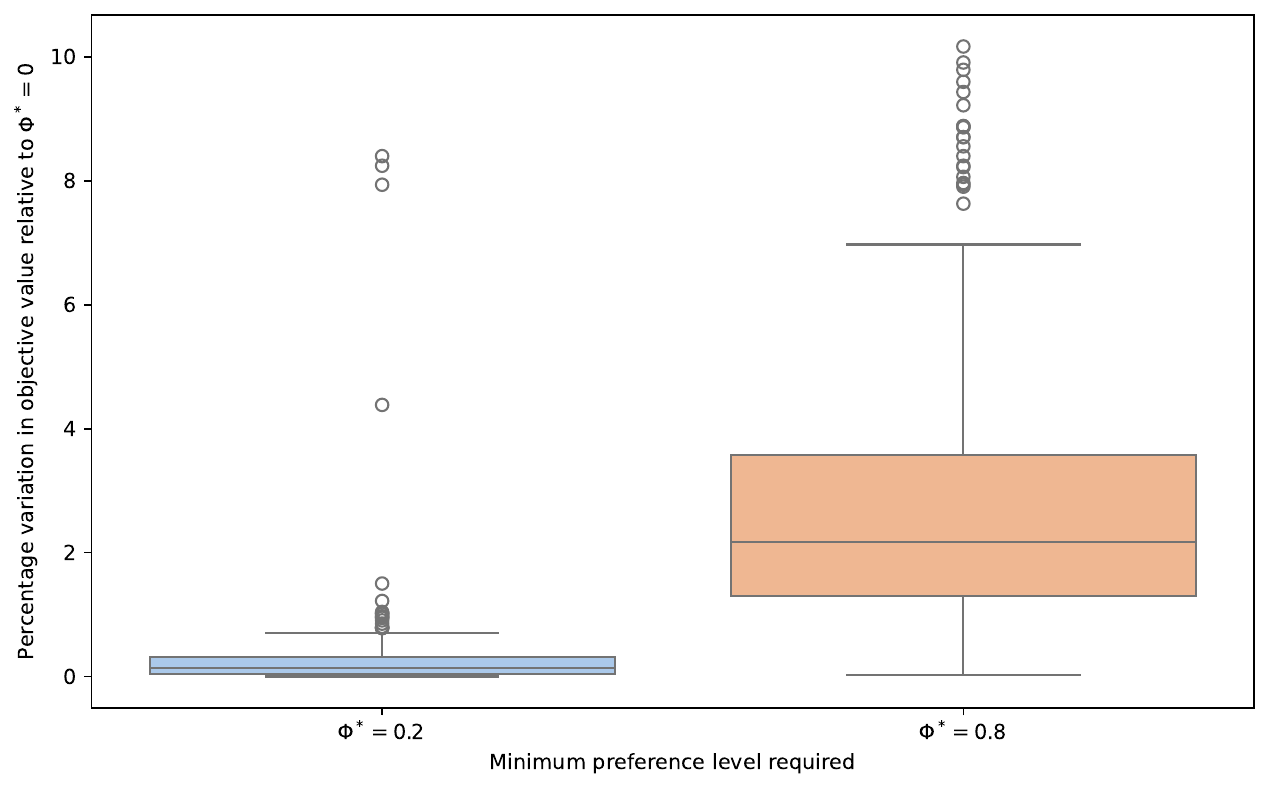}
        \caption{$n=20$.}\label{fig:boxplot_percentage20}
    \end{subfigure}~\begin{subfigure}[b]{.5\linewidth}
        \centering 
        \includegraphics[scale=0.32]{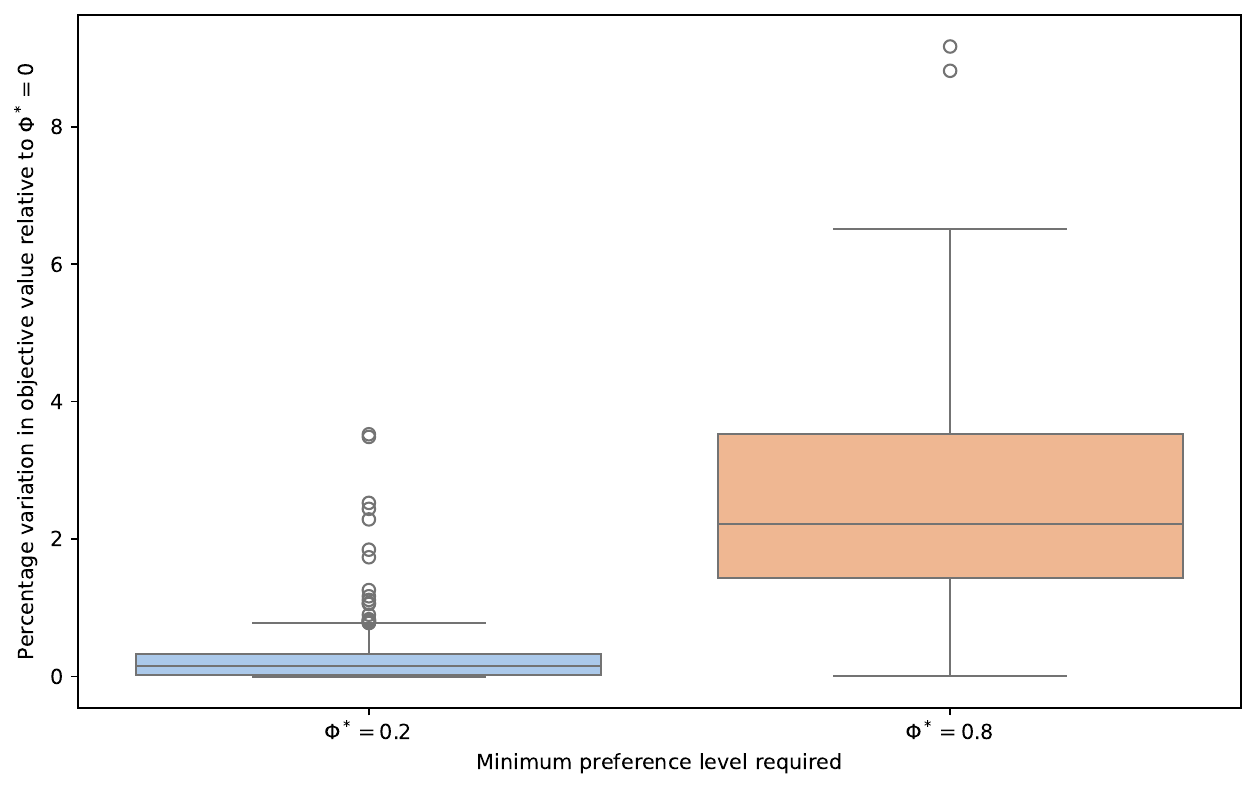}
        \caption{$n=50$.}\label{fig:boxplot_percentage50}
    \end{subfigure}\\ \begin{subfigure}[b]{.5\linewidth}
        \centering 
        \includegraphics[scale=0.32]{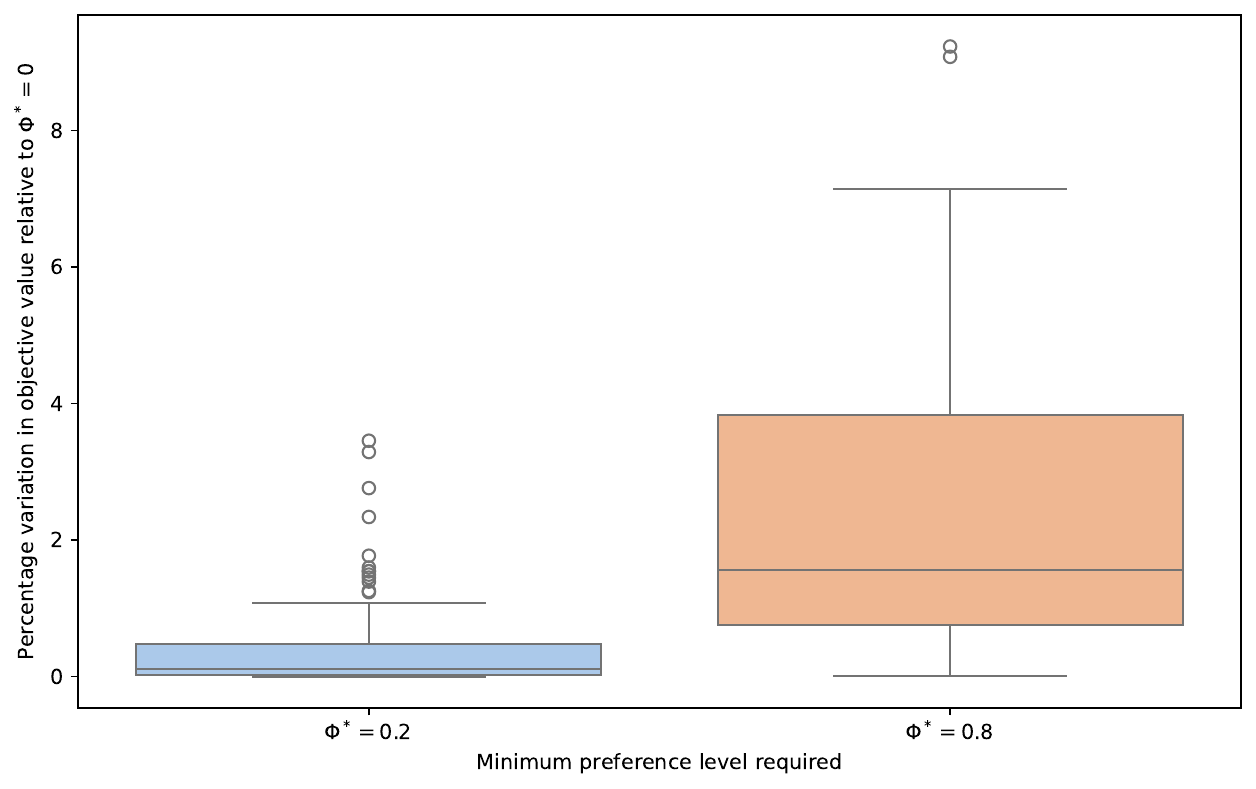}
        \caption{$n=200$.}\label{fig:boxplot_percentage200}
    \end{subfigure}~\begin{subfigure}[b]{.5\linewidth}
        \centering 
        \includegraphics[scale=0.32]{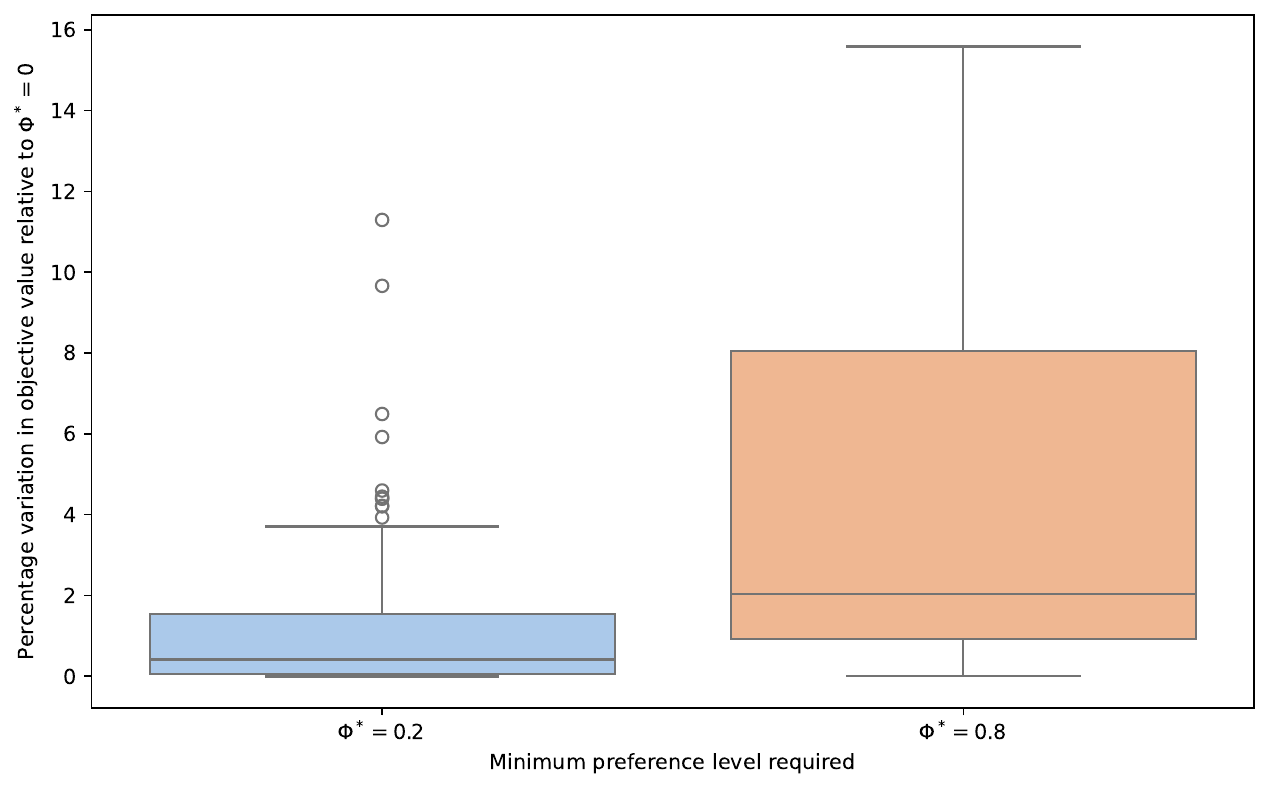}
        \caption{$n=500$.}\label{fig:boxplot_percentage500}
    \end{subfigure}~
    \caption{Boxplots of the price of efficienct of different threshold levels ($\Phi^*$) and number of regions $n$.}\label{fig:boxplot_percentage}
\end{figure}

Intrigued by this observation, in Figure \ref{fig:boxplot_percentage} we plot the \textit{Price of Efficiency}~\cite{bertsimas2011price} in the shape of boxplots representing all the instances for the two non-zero thresholds that we consider, disaggregated by the value of $n$. In our case, this value measures, for each instance, the percent deviation of the transportation cost when a given preference threshold is considered with respect to the case where no preferences are incorporated (or equivalently, when the threshold value is zero), namely:
$$
{\rm PE}(\Phi^*) = \frac{{\rm TC}(\Phi^*) - {\rm TC}(0)}{{\rm TC}(0)}, 
$$
where ${\rm TC}(\Phi^*)$ is the optimal value (transportation cost) for the instance when the threshold $\Phi^*$ is considered for the preferences ($\Phi^* \in \{0, 0.2, 0.8\}$.

As expected, the transportation cost for thresholds of 0.2 and 0.8 is positive, indicating that considering customer preferences increases this cost. Additionally, it can be observed that as the threshold increases, the cost also rises. 
The reader can also observe that the objective values (transportation costs) do not dramatically increase when incorporating preferences. The maximum relative deviation between the problems considering preferences with respect to the $\Phi^*=0$ case is at most $10\%$ for $n=20$ and $16\%$ for $n=500$. Thus, one can conclude that the incorporation of preferences, apart from satisfying the users, do not increase the operational costs assumed by the agent installing the distribution system. 

To conclude, we investigate whether changes in thresholds are statistically significant in transportation costs. With this, we aim to determine if incorporating preferences leads to statistically significant different solutions, once we know that preferences do not come with a high cost for the decision-maker.

\begin{figure}[ht]
    \begin{subfigure}[b]{.5\linewidth}
        \centering 
        \includegraphics[scale=0.48]{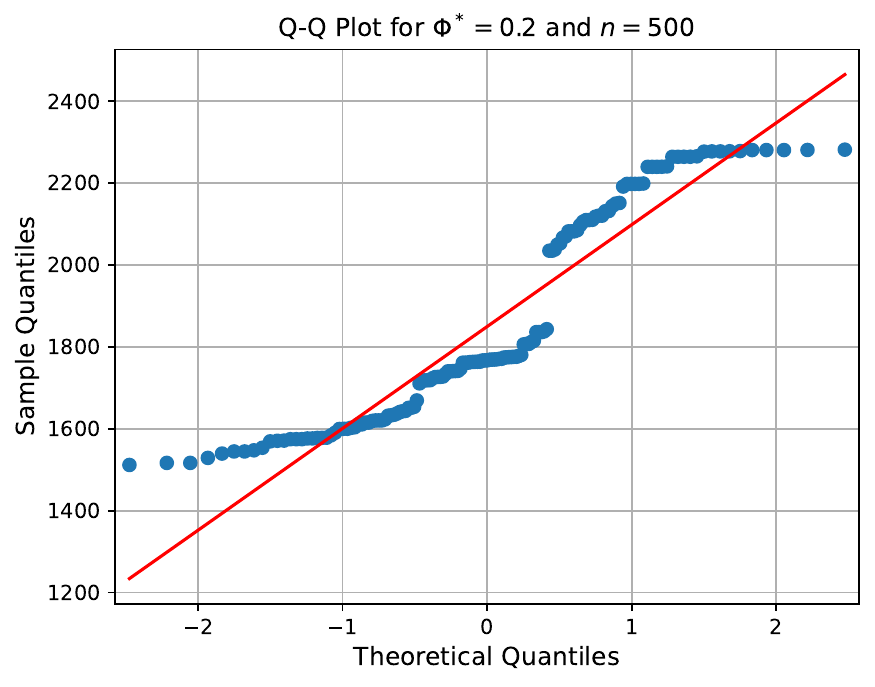}
        \caption{QQplot for $n=500$ and $\Phi^* = 0.2$.}\label{fig:normality1}
    \end{subfigure}~\begin{subfigure}[b]{.5\linewidth}
        \centering 
        \includegraphics[scale=0.4]{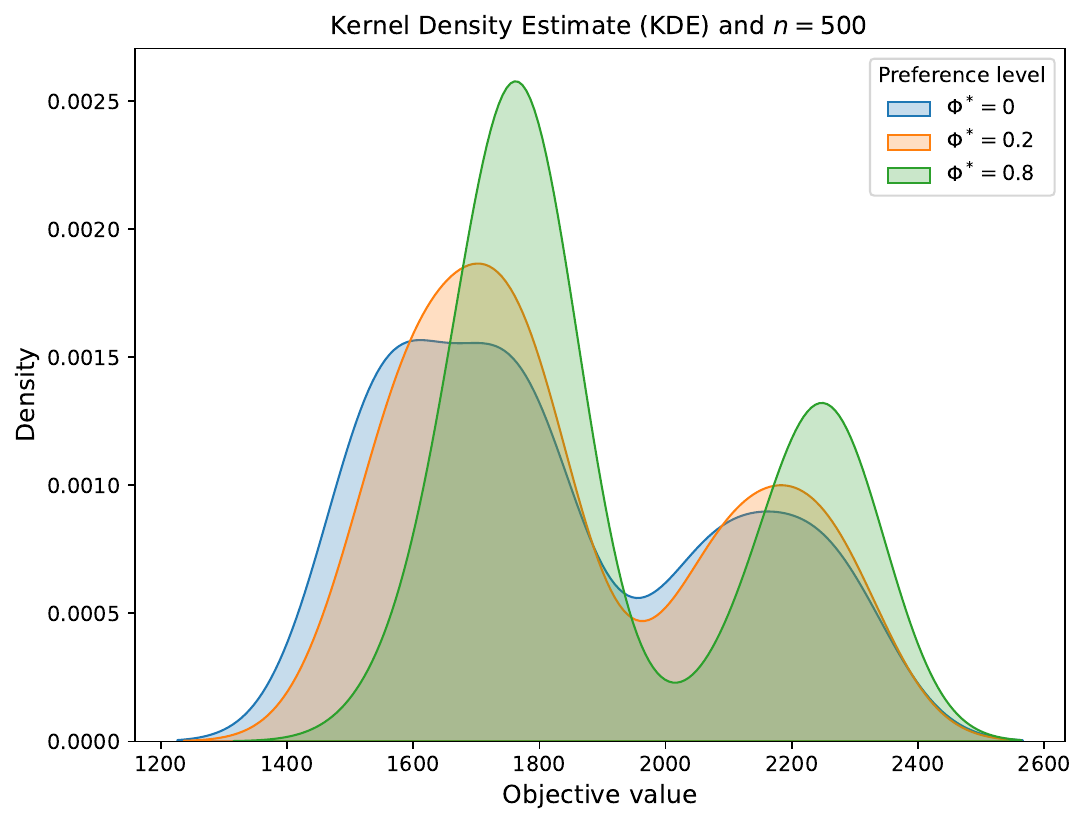}
        \caption{KDE for $n=500$.}\label{fig:normality2}
    \end{subfigure}~
    \caption{Visual demonstration of the non-normality of our data. On the left, QQ-plot between the theoretical quantiles and those of the target values for this configuration. On the right, a Kernel Density Estimate (KDE) showing the distribution for each minimum preference level.}\label{fig:normality}
\end{figure}

Since our dataset was not generated following a normal distribution, the distribution of the objective function does not follow a normal distribution in most cases. We illustrate this general trend in the  example in Figure \ref{fig:normality} with a QQ-plot that compares the quantiles of a normal distribution with the obtained values for $n=500$. This figure shows the distribution of the objective values produced by the model, concluding that our set of objective values does not follow a Gaussian distribution. 

Thus, we applied the Kruskal-Wallis test, which compares the medians among the three different thresholds. The result of this test indicates whether there are differences between the medians. 
The results of the Kruskal-Wallis test for each $n$ are summarized in Table \ref{tab:statistics}. We observe that the $p$-value decreases as the dataset size increases, with the exception of the outlier case at $n=200$. Using a significance level of $5\%$, we can conclude that there are significant differences in the medians of the objective values across the three preference levels for $n\in \{20,100,500\}$.

\begin{table}[ht]
\centering
\begin{tabular}{lrr}
    \toprule
    $n$ & Statistic & $p$-Value \\ 
    \hline
    10  & 3.004 & 0.223 \\
    20  & 6.393 & 0.041 \\
    50  & 5.486 & 0.064 \\
    100 & 10.603 & 0.005 \\
    200 & 3.884 & 0.143 \\
    500 & 20.75 & 0.0 \\
    \bottomrule
\end{tabular}
\caption{Results for the Kruskal-Wallis test on medians.} \label{tab:statistics}
\end{table}

We now perform Dunn's test analysis for these three sizes, $n \in \{20, 100, 500\}$, to determine between which pairs of sets differences exist. Table \ref{tab:dunn} presents the $p$-values from the pairwise comparisons of the means for different $\Phi^*$ combinations, for each of the three $n$ values.

In the first row, we compare $\Phi^* = 0$ with each $\Phi^* \in \{0, 0.2, 0.8\}$. We observe a $p$-value of 1 when comparing $\Phi^* = 0$ against itself, as expected.
The comparison between $\Phi^* = 0$ and $\Phi^* = 0.2$ is not  statistically significative, although the $p$-value decreases when $n=500$. Additionally, when comparing $\Phi^* = 0$ with $\Phi^* = 0.8$, the test is significant for all $n$ values, indicating that high thresholds impact the model's performance. Finally, comparing $\Phi^* = 0.2$ with $\Phi^* = 0.8$, we find that it is not significant for $n=20$, but as the dataset size increases, it becomes significant (notably at $n=500$).

\begin{table}[ht]
\centering
\adjustbox{width=\textwidth}{
\begin{tabular}{lrrrrrrrrr}
\toprule
& \multicolumn{3}{c}{$n=20$} & \multicolumn{3}{c}{$n=100$} & \multicolumn{3}{c}{$n=500$} \\
\cline{2-10}
$\Phi^*$ & 0 & 0.2 & 0.8 & 0 & 0.2 & 0.8 & 0 & 0.2 & 0.8 \\
\hline
0   & 1.0000 & 1.0000 & 0.0435 & 1.0000 & 1.0000 & 0.0047 & 1.0000 & 0.7147 & 0.0000 \\
0.2 &  & 1.0000 & 0.2243 &  & 1.0000 & 0.0717 &  & 1.0000 & 0.0039 \\
0.8 &  &  & 1.0000 &  &  & 1.0000 &  &  & 1.0000 \\
\bottomrule
\end{tabular}
}
\caption{Results of Dunn test with Bonferroni correction on medians.} \label{tab:dunn}
\end{table}

Therefore, from this statistical analysis, we draw two interesting conclusions regarding the applicability of our model to real-world cases:
\begin{enumerate}
\item The incorporation of preferences allows to improve the customer's satisfaction, by a small increase in the operational costs. In our extensive  computational experience, the costs increase, in average,  $4-5\%$, with a maximum increase of $16\%$  in a few cases. Thus, the incorporation of preferences is justified from an economic perspective.
\item From a practical perspective, this incorporation has a better performance when setting high thresholds for small-sized datasets. As the sample size increases, it becomes more practical to apply different threshold values for locational decisions, allowing for more flexibility for real-world  applications.
\end{enumerate}

\section{Conclusions and Future Research}\label{sec:conclusions}

This work analyzes the incorporation of preferences within demand regions for continuous facility location problems. We analyze mathematical optimization models for two interesting cases: the general case where collocation is not considered, and the collocation framework where locating services at the same position might be beneficial for overlapping regions. The models are rewritten as MISOCO problems that can be solved using off-the-shelf software. Additionally, we provide up to five different preference functions that can be applied to real-world cases, together with the details about how to incorporate them into the mathematical optimization models.

The paper presents two main lines of work: computational and qualitative. In the computational section, we achieved optimal solutions for problems with up to 500 demand points in the single-facility case. However, the problem becomes computationally costly for the multi-facility case, even for medium-sized instances.

On the qualitative side, we address whether incorporating preferences adds value to the decision-making process. The answer is a resounding yes, but primarily for large datasets. Specifically, when dealing with small sets of demand regions, the results indicate no significant differences between incorporating preferences or not. However, as $n$ increases, we observe significant differences.

The work done encourages some other research lines. First, in this work, we formulated the problem assuming a general benefit of the whole demand regions; which makes sense in the case in which the decision maker of the locations of all service points and facilities is the same. However, there may exist different decision-makers deciding, separately and hierarchically, about the locations. Such a case could be addressed by using multilevel optimization. This will be the topic of a forthcoming paper.

Another possible extension of this work is by finding, apart from the facility and serving points locations, the best routes between the serving points and the facilities; instead of assuming a direct connection between each serving point and the facility. That is, incorporating a Traveling Salesman Problem or Vehicle Routing Problem in the Continuous Facility Location Problem with Demand Regions and Preferences. This extension will also be studied in a future paper.

\section*{Acknowledgements}

Authors thank the ``Servicio de Supercomputacion" from Universidad de Granada \sloppy (\url{https://supercomputacion.ugr.es}) for providing computing time on the albaicin supercomputer.
This research has been partially supported by grant PID2020-114594GB-C21 funded by MICIU/AEI/ 10.13039/501100011033, grant RED2022-134149-T funded by MICIU/AEI /10.13039/501100011033(Thematic Network on Location Science and Related Problems), and the IMAG-María de Maeztu grant
CEX2020-001105-M/AEI/10.13039/501100011033.

\bibliographystyle{elsarticle-harv} 
\bibliography{00_references}

\begin{thebibliography}{38}
\expandafter\ifx\csname natexlab\endcsname\relax\def\natexlab#1{#1}\fi
\providecommand{\url}[1]{\texttt{#1}}
\providecommand{\href}[2]{#2}
\providecommand{\path}[1]{#1}
\providecommand{\DOIprefix}{doi:}
\providecommand{\ArXivprefix}{arXiv:}
\providecommand{\URLprefix}{URL: }
\providecommand{\Pubmedprefix}{pmid:}
\providecommand{\doi}[1]{\href{http://dx.doi.org/#1}{\path{#1}}}
\providecommand{\Pubmed}[1]{\href{pmid:#1}{\path{#1}}}
\providecommand{\bibinfo}[2]{#2}
\ifx\xfnm\relax \def\xfnm[#1]{\unskip,\space#1}\fi
\bibitem[{Bertsimas et~al.(2011)Bertsimas, Farias and
  Trichakis}]{bertsimas2011price}
\bibinfo{author}{Bertsimas, D.}, \bibinfo{author}{Farias, V.F.},
  \bibinfo{author}{Trichakis, N.}, \bibinfo{year}{2011}.
\newblock \bibinfo{title}{The price of fairness}.
\newblock \bibinfo{journal}{Operations research} \bibinfo{volume}{59},
  \bibinfo{pages}{17--31}.
\bibitem[{Blanco(2019)}]{blanco2019ordered}
\bibinfo{author}{Blanco, V.}, \bibinfo{year}{2019}.
\newblock \bibinfo{title}{Ordered p-median problems with neighbourhoods}.
\newblock \bibinfo{journal}{Computational Optimization and Applications}
  \bibinfo{volume}{73}, \bibinfo{pages}{603--645}.
\bibitem[{Blanco et~al.(2017a)Blanco, Fern{\'a}ndez and
  Puerto}]{blanco2017minimum}
\bibinfo{author}{Blanco, V.}, \bibinfo{author}{Fern{\'a}ndez, E.},
  \bibinfo{author}{Puerto, J.}, \bibinfo{year}{2017}a.
\newblock \bibinfo{title}{Minimum spanning trees with neighborhoods:
  Mathematical programming formulations and solution methods}.
\newblock \bibinfo{journal}{European Journal of Operational Research}
  \bibinfo{volume}{262}, \bibinfo{pages}{863--878}.
\bibitem[{Blanco and G{\'a}zquez(2023)}]{blanco2023fairness}
\bibinfo{author}{Blanco, V.}, \bibinfo{author}{G{\'a}zquez, R.},
  \bibinfo{year}{2023}.
\newblock \bibinfo{title}{Fairness in maximal covering location problems}.
\newblock \bibinfo{journal}{Computers \& Operations Research}
  \bibinfo{volume}{157}, \bibinfo{pages}{106287}.
\bibitem[{Blanco et~al.(2023)Blanco, G{\'a}zquez and Saldanha-da
  Gama}]{blanco2023multi}
\bibinfo{author}{Blanco, V.}, \bibinfo{author}{G{\'a}zquez, R.},
  \bibinfo{author}{Saldanha-da Gama, F.}, \bibinfo{year}{2023}.
\newblock \bibinfo{title}{Multi-type maximal covering location problems:
  Hybridizing discrete and continuous problems}.
\newblock \bibinfo{journal}{European journal of operational research}
  \bibinfo{volume}{307}, \bibinfo{pages}{1040--1054}.
\bibitem[{Blanco et~al.(2024)Blanco, Mar{\'\i}n and Puerto}]{blanco2024intra}
\bibinfo{author}{Blanco, V.}, \bibinfo{author}{Mar{\'\i}n, A.},
  \bibinfo{author}{Puerto, J.}, \bibinfo{year}{2024}.
\newblock \bibinfo{title}{Intra-facility equity in discrete and continuous
  p-facility location problems}.
\newblock \bibinfo{journal}{Computers \& Operations Research}
  \bibinfo{volume}{162}, \bibinfo{pages}{106487}.
\bibitem[{Blanco and Mart{\'\i}nez-Ant{\'o}n(2024)}]{blanco2023minimal}
\bibinfo{author}{Blanco, V.}, \bibinfo{author}{Mart{\'\i}nez-Ant{\'o}n, M.},
  \bibinfo{year}{2024}.
\newblock \bibinfo{title}{On minimal extended representations of generalized
  power cones}.
\newblock \bibinfo{journal}{SIAM J. Optimization} \bibinfo{volume}{To appear}.
\bibitem[{Blanco and Puerto(2022)}]{blanco2022hub}
\bibinfo{author}{Blanco, V.}, \bibinfo{author}{Puerto, J.},
  \bibinfo{year}{2022}.
\newblock \bibinfo{title}{On hub location problems in geographically flexible
  networks}.
\newblock \bibinfo{journal}{International Transactions in Operational Research}
  \bibinfo{volume}{29}, \bibinfo{pages}{2226--2249}.
\bibitem[{Blanco et~al.(2016)Blanco, Puerto and Ben-Ali}]{BEP16}
\bibinfo{author}{Blanco, V.}, \bibinfo{author}{Puerto, J.},
  \bibinfo{author}{Ben-Ali, S.E.H.}, \bibinfo{year}{2016}.
\newblock \bibinfo{title}{Continuous multifacility ordered median location
  problems}.
\newblock \bibinfo{journal}{European Journal of Operational Research}
  \bibinfo{volume}{250}, \bibinfo{pages}{56--64}.
\bibitem[{Blanco et~al.(2014)Blanco, Puerto and El~Haj Ben~Ali}]{BEP14}
\bibinfo{author}{Blanco, V.}, \bibinfo{author}{Puerto, J.},
  \bibinfo{author}{El~Haj Ben~Ali, S.}, \bibinfo{year}{2014}.
\newblock \bibinfo{title}{Revisiting several problems and algorithms in
  continuous location with $\ell_\tau$-norms}.
\newblock \bibinfo{journal}{Computational Optimization and Applications}
  \bibinfo{volume}{58}, \bibinfo{pages}{563--595}.
\bibitem[{Blanco et~al.(2017b)Blanco, Puerto and Ponce}]{blanco2017continuous}
\bibinfo{author}{Blanco, V.}, \bibinfo{author}{Puerto, J.},
  \bibinfo{author}{Ponce, D.}, \bibinfo{year}{2017}b.
\newblock \bibinfo{title}{Continuous location under the effect of
  ‘refraction’}.
\newblock \bibinfo{journal}{Mathematical Programming} \bibinfo{volume}{161},
  \bibinfo{pages}{33--72}.
\bibitem[{Camacho-Vallejo et~al.(2014)Camacho-Vallejo, Casas-Ram{\'\i}rez and
  Miranda}]{camacho2014p}
\bibinfo{author}{Camacho-Vallejo, J.F.}, \bibinfo{author}{Casas-Ram{\'\i}rez,
  M.}, \bibinfo{author}{Miranda, P.}, \bibinfo{year}{2014}.
\newblock \bibinfo{title}{The p-median bilevel problem under preferences of the
  customers}.
\newblock \bibinfo{journal}{Recent Advances in Theory, Methods and Practice of
  Operations Research} , \bibinfo{pages}{121--127}.
\bibitem[{Carrizosa et~al.(1998)Carrizosa, Mu{\~n}oz-M{\'a}rquez and
  Puerto}]{carrizosa1998weber}
\bibinfo{author}{Carrizosa, E.}, \bibinfo{author}{Mu{\~n}oz-M{\'a}rquez, M.},
  \bibinfo{author}{Puerto, J.}, \bibinfo{year}{1998}.
\newblock \bibinfo{title}{The weber problem with regional demand}.
\newblock \bibinfo{journal}{European Journal of Operational Research}
  \bibinfo{volume}{104}, \bibinfo{pages}{358--365}.
\bibitem[{Casas-Ram{\'\i}rez and Camacho-Vallejo(2017)}]{casas2017solving}
\bibinfo{author}{Casas-Ram{\'\i}rez, M.S.}, \bibinfo{author}{Camacho-Vallejo,
  J.F.}, \bibinfo{year}{2017}.
\newblock \bibinfo{title}{Solving the p-median bilevel problem with order
  through a hybrid heuristic}.
\newblock \bibinfo{journal}{Applied Soft Computing} \bibinfo{volume}{60},
  \bibinfo{pages}{73--86}.
\bibitem[{Chanta et~al.(2014)Chanta, Mayorga and McLay}]{chanta2014minimum}
\bibinfo{author}{Chanta, S.}, \bibinfo{author}{Mayorga, M.E.},
  \bibinfo{author}{McLay, L.A.}, \bibinfo{year}{2014}.
\newblock \bibinfo{title}{The minimum p-envy location problem with requirement
  on minimum survival rate}.
\newblock \bibinfo{journal}{Computers \& Industrial Engineering}
  \bibinfo{volume}{74}, \bibinfo{pages}{228--239}.
\bibitem[{Cobb and Douglas(1928)}]{cobb1928theory}
\bibinfo{author}{Cobb, C.W.}, \bibinfo{author}{Douglas, P.H.},
  \bibinfo{year}{1928}.
\newblock \bibinfo{title}{A theory of production}.
\newblock \bibinfo{journal}{American Economic Review} \bibinfo{volume}{18},
  \bibinfo{pages}{139--165}.
\bibitem[{Disser et~al.(2014)Disser, Mihal{\'a}k, Montanari and
  Widmayer}]{disser2014rectilinear}
\bibinfo{author}{Disser, Y.}, \bibinfo{author}{Mihal{\'a}k, M.},
  \bibinfo{author}{Montanari, S.}, \bibinfo{author}{Widmayer, P.},
  \bibinfo{year}{2014}.
\newblock \bibinfo{title}{Rectilinear shortest path and rectilinear minimum
  spanning tree with neighborhoods}, in: \bibinfo{booktitle}{International
  Symposium on Combinatorial Optimization}, \bibinfo{organization}{Springer}.
  pp. \bibinfo{pages}{208--220}.
\bibitem[{Drezner(2022)}]{drezner2022continuous}
\bibinfo{author}{Drezner, Z.}, \bibinfo{year}{2022}.
\newblock \bibinfo{title}{Continuous facility location problems}, in:
  \bibinfo{booktitle}{The Palgrave handbook of operations research}.
  \bibinfo{publisher}{Springer}, pp. \bibinfo{pages}{269--306}.
\bibitem[{Drezner and Weslowsky(1980)}]{drezner1980optimal}
\bibinfo{author}{Drezner, Z.}, \bibinfo{author}{Weslowsky, G.},
  \bibinfo{year}{1980}.
\newblock \bibinfo{title}{Optimal location of a facility relative to area
  demands}.
\newblock \bibinfo{journal}{Naval Research Logistics Quarterly}
  \bibinfo{volume}{27}, \bibinfo{pages}{199--206}.
\bibitem[{Espejo et~al.(2009)Espejo, Mar{\'\i}n, Puerto and
  Rodr{\'\i}guez-Ch{\'\i}a}]{espejo2009comparison}
\bibinfo{author}{Espejo, I.}, \bibinfo{author}{Mar{\'\i}n, A.},
  \bibinfo{author}{Puerto, J.}, \bibinfo{author}{Rodr{\'\i}guez-Ch{\'\i}a,
  A.M.}, \bibinfo{year}{2009}.
\newblock \bibinfo{title}{A comparison of formulations and solution methods for
  the minimum-envy location problem}.
\newblock \bibinfo{journal}{Computers \& Operations Research}
  \bibinfo{volume}{36}, \bibinfo{pages}{1966--1981}.
\bibitem[{Espejo et~al.(2022)Espejo, P{\'a}ez, Puerto and
  Rodr{\'\i}guez-Ch{\'\i}a}]{espejo2022minimum}
\bibinfo{author}{Espejo, I.}, \bibinfo{author}{P{\'a}ez, R.},
  \bibinfo{author}{Puerto, J.}, \bibinfo{author}{Rodr{\'\i}guez-Ch{\'\i}a, A.},
  \bibinfo{year}{2022}.
\newblock \bibinfo{title}{Minimum cost b-matching problems with neighborhoods}.
\newblock \bibinfo{journal}{Computational Optimization and Applications}
  \bibinfo{volume}{83}, \bibinfo{pages}{525--553}.
\bibitem[{Espejo et~al.(2023)Espejo, P{\'a}ez, Puerto and
  Rodr{\'\i}guez-Ch{\'\i}a}]{espejo2023facility}
\bibinfo{author}{Espejo, I.}, \bibinfo{author}{P{\'a}ez, R.},
  \bibinfo{author}{Puerto, J.}, \bibinfo{author}{Rodr{\'\i}guez-Ch{\'\i}a, A.},
  \bibinfo{year}{2023}.
\newblock \bibinfo{title}{Facility location problems on graphs with non-convex
  neighborhoods}.
\newblock \bibinfo{journal}{Computers \& Operations Research}
  \bibinfo{volume}{159}, \bibinfo{pages}{106356}.
\bibitem[{Grossmann(2002)}]{grossmann2002review}
\bibinfo{author}{Grossmann, I.E.}, \bibinfo{year}{2002}.
\newblock \bibinfo{title}{Review of nonlinear mixed-integer and disjunctive
  programming techniques}.
\newblock \bibinfo{journal}{Optimization and engineering} \bibinfo{volume}{3},
  \bibinfo{pages}{227--252}.
\bibitem[{Hakimi(1964)}]{hakimi1964optimum}
\bibinfo{author}{Hakimi, S.L.}, \bibinfo{year}{1964}.
\newblock \bibinfo{title}{Optimum locations of switching centers and the
  absolute centers and medians of a graph}.
\newblock \bibinfo{journal}{Operations research} \bibinfo{volume}{12},
  \bibinfo{pages}{450--459}.
\bibitem[{Hakimi(1983)}]{hakimi1983locating}
\bibinfo{author}{Hakimi, S.L.}, \bibinfo{year}{1983}.
\newblock \bibinfo{title}{On locating new facilities in a competitive
  environment}.
\newblock \bibinfo{journal}{European Journal of Operational Research}
  \bibinfo{volume}{12}, \bibinfo{pages}{29--35}.
\bibitem[{Hanjoul and Peeters(1987)}]{hanjoul1987facility}
\bibinfo{author}{Hanjoul, P.}, \bibinfo{author}{Peeters, D.},
  \bibinfo{year}{1987}.
\newblock \bibinfo{title}{A facility location problem with clients' preference
  orderings}.
\newblock \bibinfo{journal}{Regional Science and Urban Economics}
  \bibinfo{volume}{17}, \bibinfo{pages}{451--473}.
\bibitem[{Hansen et~al.(1987)Hansen, Labb{\'e}, Peeters and
  Thisse}]{hansen1987single}
\bibinfo{author}{Hansen, P.}, \bibinfo{author}{Labb{\'e}, M.},
  \bibinfo{author}{Peeters, D.}, \bibinfo{author}{Thisse, J.F.},
  \bibinfo{year}{1987}.
\newblock \bibinfo{title}{Single facility location on networks}, in:
  \bibinfo{booktitle}{North-Holland mathematics studies}.
  \bibinfo{publisher}{Elsevier}. volume \bibinfo{volume}{132}, pp.
  \bibinfo{pages}{113--145}.
\bibitem[{K{\"u}{\c{c}}{\"u}kayd{\i}n and Aras(2020)}]{kuccukaydin2020gradual}
\bibinfo{author}{K{\"u}{\c{c}}{\"u}kayd{\i}n, H.}, \bibinfo{author}{Aras, N.},
  \bibinfo{year}{2020}.
\newblock \bibinfo{title}{Gradual covering location problem with multi-type
  facilities considering customer preferences}.
\newblock \bibinfo{journal}{Computers \& Industrial Engineering}
  \bibinfo{volume}{147}, \bibinfo{pages}{106577}.
\bibitem[{Laporte et~al.(2019)Laporte, Nickel and
  {Saldanha-da-Gama}}]{LaporteNickelSaldanha-da-Gama:2019}
\bibinfo{author}{Laporte, G.}, \bibinfo{author}{Nickel, S.},
  \bibinfo{author}{{Saldanha-da-Gama}, F.}, \bibinfo{year}{2019}.
\newblock \bibinfo{title}{Location Science}.
\newblock \bibinfo{publisher}{Springer International Publishing, 2nd edition}.
\bibitem[{Leontief(1941)}]{leontief1941}
\bibinfo{author}{Leontief, W.}, \bibinfo{year}{1941}.
\newblock \bibinfo{title}{The structure of the american economy, 1919-1929: An
  empirical application of equilibrium analysis}.
\newblock \bibinfo{journal}{Harvard University Press} .
\bibitem[{Megiddo and Supowit(1984)}]{megiddo1984complexity}
\bibinfo{author}{Megiddo, N.}, \bibinfo{author}{Supowit, K.J.},
  \bibinfo{year}{1984}.
\newblock \bibinfo{title}{On the complexity of some common geometric location
  problems}.
\newblock \bibinfo{journal}{SIAM journal on computing} \bibinfo{volume}{13},
  \bibinfo{pages}{182--196}.
\bibitem[{Michelot(1987)}]{michelot1987localization}
\bibinfo{author}{Michelot, C.}, \bibinfo{year}{1987}.
\newblock \bibinfo{title}{Localization in multifacility location theory}.
\newblock \bibinfo{journal}{European Journal of Operational Research}
  \bibinfo{volume}{31}, \bibinfo{pages}{177--184}.
\bibitem[{Puerto et~al.(2018)Puerto, Ricca and Scozzari}]{puerto2018extensive}
\bibinfo{author}{Puerto, J.}, \bibinfo{author}{Ricca, F.},
  \bibinfo{author}{Scozzari, A.}, \bibinfo{year}{2018}.
\newblock \bibinfo{title}{Extensive facility location problems on networks: an
  updated review}.
\newblock \bibinfo{journal}{Top} \bibinfo{volume}{26},
  \bibinfo{pages}{187--226}.
\bibitem[{Schoenwitz et~al.(2017)Schoenwitz, Potter, Gosling and
  Naim}]{schoenwitz2017product}
\bibinfo{author}{Schoenwitz, M.}, \bibinfo{author}{Potter, A.},
  \bibinfo{author}{Gosling, J.}, \bibinfo{author}{Naim, M.},
  \bibinfo{year}{2017}.
\newblock \bibinfo{title}{Product, process and customer preference alignment in
  prefabricated house building}.
\newblock \bibinfo{journal}{International Journal of Production Economics}
  \bibinfo{volume}{183}, \bibinfo{pages}{79--90}.
\bibitem[{Solow(1956)}]{solow}
\bibinfo{author}{Solow, R.M.}, \bibinfo{year}{1956}.
\newblock \bibinfo{title}{A contribution to the theory of economic growth}.
\newblock \bibinfo{journal}{The quarterly journal of economics}
  \bibinfo{volume}{70}, \bibinfo{pages}{65--94}.
\bibitem[{Wang et~al.(2024)Wang, Wong, Shi and Yuen}]{wang2024investigation}
\bibinfo{author}{Wang, X.}, \bibinfo{author}{Wong, Y.D.}, \bibinfo{author}{Shi,
  W.}, \bibinfo{author}{Yuen, K.F.}, \bibinfo{year}{2024}.
\newblock \bibinfo{title}{An investigation on consumers' preferences for parcel
  deliveries: applying consumer logistics in omni-channel shopping}.
\newblock \bibinfo{journal}{The International Journal of Logistics Management}
  \bibinfo{volume}{35}, \bibinfo{pages}{557--576}.
\bibitem[{Weber(1909)}]{weber1909standort}
\bibinfo{author}{Weber, A.}, \bibinfo{year}{1909}.
\newblock \bibinfo{title}{{\"U}ber den standort der industrien. 1. teil: Reine
  theorie des standorts}.
\newblock \bibinfo{journal}{JCB Mohr, T{\"u}bingen,(English ed. by CJ
  Friedrichs, Univ. Chicago Press, 1929)} .
\bibitem[{Yura(1994)}]{yura1994production}
\bibinfo{author}{Yura, K.}, \bibinfo{year}{1994}.
\newblock \bibinfo{title}{Production scheduling to satisfy worker's preferences
  for days off and overtime under due-date constraints}.
\newblock \bibinfo{journal}{International Journal of Production Economics}
  \bibinfo{volume}{33}, \bibinfo{pages}{265--270}.

\end{thebibliography}

\end{document}